\documentclass{amsart}
\usepackage{psfig}
\usepackage{amssymb, amsmath,amscd}
\def\epsffile#1{\psfig{file=#1,silent=}}

\newcommand{\R}{\mathbb R}
\newcommand{\Z}{\mathbb Z}
\def\Ga{\alpha}

\def\GD{\Delta}
\def\Ge{\varepsilon}
\def\Gf{\phi}

\def\Gg{\gamma}
\def\GG{\Gamma}
\def\Gl{\lambda}
\def\Go{\omega}

\def\Gp{\pi}
\def\Gs{\sigma}
\def\GS{\Sigma}

\def\Gz{\zeta}
\def\sign{\operatorname{sign}}
\def\p{\partial}
\theoremstyle{plain}
\newtheorem{thm}{Theorem}[section]

\newtheorem{lem}{Lemma}
\newtheorem{cor}[thm]{Corollary}
\theoremstyle{definition}
\newtheorem{ex}[thm]{Example}

\theoremstyle{remark}
\newtheorem{rem}[thm]{Remark}
\newcommand{\abz}{\newline\indent}
\newcommand{\cal}{\mathcal}

\newcommand{\sminus}{\smallsetminus}
\def\id{\operatorname{id}}
\begin{document}
\title{On the Casson Knot Invariant}
\author[M. Polyak]{Michael Polyak}
\address{M.~Polyak,
School of Mathematics, Tel-Aviv University,
69978 Tel-Aviv, Israel}
\email{polyak@math.tau.ac.il}
\author[O. Viro]{Oleg Viro}
\address{O.~Viro, Department of Mathematics, Uppsala University S-751 06
Uppsala, Sweden;\abz
POMI, Fontanka 27, St.Petersburg, 191011, Russia.}
\email{oleg@math.uu.se} 
\maketitle


In our previous paper \cite{PV} we introduced a new type  of combinatorial
formulas for Vassiliev knot invariants and presented lots of formulas
of this type. To the best of our knowledge, these formulas are by far the
simplest and the most practical for computational purposes.  Since then
Goussarov  has proved the main conjecture
formulated in \cite{PV}:  any Vassiliev knot invariant can be described
by such a formula, see \cite{GPV}.

In \cite{PV} the examples of formulas were presented in a formal way,
without proofs or even explanations of the ideas. We promised to
interpret the invariants as degrees of some maps in a forthcoming paper
and mentioned that it was this viewpoint that motivated the whole our
investigations and appeared to be a rich source of various special
formulas.

In a sense, this viewpoint was not new. Quite the contrary, this is the
most classical way to think on knot invariants. Indeed, a
classical definition of a knot invariant runs as follows: some
geometric construction gives an auxiliary space and then the machinery
of algebraic topology is applied to this space to produce a number (or
a quadratic form, a group, etc.). This scheme was almost forgotten in 
the eighties, when quantum invariants appeared. An auxiliary space was
replaced by a combinatorial object (like knot diagram or closed braid
presentation of a knot), while algebraic topology was replaced by
representation theory and statistical mechanics. Vassiliev invariants
and calculation of the quantum invariants in terms of Vassiliev
invariants recovered the role of algebraic topology, but it is applied
to the space of all knots, rather than to a space manufactured from a
single knot. Presentations of Vassiliev invariants as degrees of maps
would completely rehabilitate the classical approach.

However, this is not our main intention. Various presentations for
Vassiliev invariants reveal a rich geometric contents. The usual
benefits of presenting some quantity as a degree of a map are that
a degree is easy to calculate by various methods and, furthermore, 
degrees are manifestly invariant under various kinds of deformations.

We were primarily motivated by the well-known case of the linking 
number. It is the simplest Vassiliev invariant of links. 
The linking number can be computed in many different ways, see e.g 
\cite{Rolfsen}. However, all formulas can be obtained from a single 
one: the linking number of a pair of circles is the degree of the map 
of a configuration space of pairs (a point on one circle, a point on 
the other circle) to $S^2$ defined by 
$(x,y)\mapsto \frac{x-y}{|x-y|}$. 
Both the Gauss integral formula, and the combinatorial formulas in 
terms of a diagram are deduced from this interpretation via various 
methods for calculation of a degree.

We discovered that this situation is reproduced in the case of
Vassiliev invariants of higher degree. 
Both integral formulas found by Kontsevich \cite{Konts} and Bar-Natan 
\cite{Bar-Nat} (of the Knizhnik-Zamolodchikov and the Chern-Simons 
type), and the combinatorial formulas of Lannes \cite{Lannes} and our 
note \cite{PV}, can be deduced from a presentation of an invariant as 
the degree of a similar map.  
However, since the configuration spaces are getting more complicated 
as the degree increases, the number of various formulas is getting 
surprisingly large. In fact our approach is close to the one of 
Bott-Taubes \cite{BT}, however we do not restrict ourselves to
integral formulas of the Chern-Simons type, but rather try to include
formulas of all types in this scheme.

In this paper we focus on the simplest Vassiliev knot invariant $v_2$.
This invariant is of degree 2. It can be characterized as the unique 
Vassiliev invariant of degree 2 which takes values $0$ on the unknot 
and $1$ on a trefoil. It was known, however, long time before this 
characterization became possible (when Vassiliev invariants were 
introduced).
Indeed, it can be defined as $\frac12\GD_K''(1)$, the half of the 
value at 1 of the second derivative of the Alexander polynomial (or 
as the coefficient of the quadratic term of the Conway polynomial). 
It is the invariant which plays a key role in the surgery
formula for the Casson invariant of homology spheres \cite{AM}.
Following a recent folklore tradition, we shall call this knot
invariant also the {\it Casson invariant\/} or the {\it Casson knot
invariant\/}, when there is a danger of confusion with the Casson
invariant of homology spheres.

We decided to devote this paper to the Casson knot invariant for
several reasons. This is the simplest knot invariant of finite
type. On the other hand, it is related to many phenomena.
For instance, its reduction modulo 2 is the Arf invariant, which is 
the only invariant of finite degree which is invariant under knot 
cobordisms \cite{Ng}.
Furthermore, the Casson knot invariant is related to Arnold's 
invariants of generic plane curves \cite{Polyak} and \cite{LW}.
It appears as well in the theory of Casson invariant for homology
spheres and plays an important role in the recent progress on finite
degree invariants of 3-manifolds.
Technically, all phenomena and problems connected to an 
interpretation of Vassiliev knot invariants as degrees of maps 
arise already in the case of Casson knot invariant.
Recently it became clear how to treat the general case.
One can consider the universal invariant taking values in the algebra
$\overrightarrow{\cal G}$ introduced in \cite{P?} and based on
constructions of the present paper with orientations of configuration
spaces defined as in S.\ Poirier \cite{Poi}.

We postpone a presentation of the universal invariant as a generalized
degree of a map to a forthcoming paper \cite{P?}. Here we concentrate
on the geometry related to the Casson invariant. This allows us to
consider all remarkable geometric constructions and phenomena which would be
inevitably omitted in any paper dedicated to a construction of
the universal invariant.

We begin with our combinatorial formula announced in \cite{PV}. 
It is proved according to a traditional combinatorial scheme on 
the base of general definition of Vassiliev invariants. 
In fact the main ingredient of this proof is Kauffman's calculation 
\cite{Kauff} of the second coefficient of the Conway polynomial. 
We use the formula to prove an upper bound on the Casson knot 
invariant via the  number of double points of a knot diagram.

This is done in  Section \ref{s1}. Then we proceed to the main subject
and construct the configuration spaces and their maps. The first
attempt in Section \ref{s2} gives rise to an interpretation of the
Casson invariant as a local degree. The disadvantage of this
interpretation is that it is restricted to the case when the knot is 
in a general position with respect to a fixed direction of
projection, and hence it is not manifestly invariant under isotopy 
and does not lead immediately to other combinatorial formulas.

In Section \ref{s3} we enlarge the source space by adding several
patches which are similar to the configuration spaces appearing in
Chern-Simons theory. Although the new space still has a boundary, the
boundary is mapped to a fixed hypersurface of the target manifold.
Thus $v_2$ gets an interpretation as a global (though relative) degree
of a map.

In Section \ref{s4} we derive new combinatorial formulas for $v_2$
taking other regular values of the map constructed in Section \ref{s3}.
In particular, this leads to a calculation in terms of associators
which appear in a presentation of the knot diagram as a nonassociative
tangle.

In Section \ref{s5} we discuss other configuration spaces and
presentations of the Casson knot invariant as the degree of the
corresponding maps.
Various methods to compute the corresponding degrees are used   to
derive new combinatorial and integral formulas.

An essential part of this work was done when the first author was 
visiting the Max-Planck-Institut f\"ur Mathematik in Bonn, which he
wishes to thank for its hospitality.

\section{A Gauss Diagram Formula for the Casson Invariant}\label{s1}

\subsection{Gauss diagrams}\label{s1.1} 
\begin{figure}[h]
\centerline{\epsffile{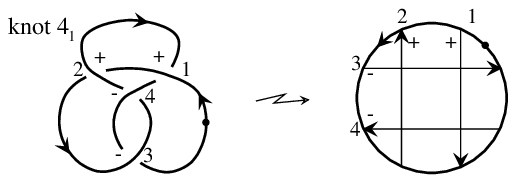}}
\caption{}
\label{f2}
\end{figure}

A knot diagram is a
generic immersion of circle to plane, enhanced by the information 
on overpasses and underpasses at double points.
A generic immersion of a circle to plane is characterized by its
{\it Gauss diagram.\/} The Gauss diagram is the immersing circle with
the preimages of each double point connected with a chord.
To incorporate the information on overpasses and underpasses, we 
orient each chord from the upper branch to the lower branch. 
Furthermore, each chord $c$ is equipped with the sign $\Ge(c)$ of 
the corresponding  double point (local writhe number). 
See Figure \ref{f2}. 
We call the result a {\it Gauss diagram \/} of the knot.

By a {\it based\/} Gauss diagram we mean a Gauss diagram with a marked
point on the circle, distinct from the end points of the chords.

\subsection{The Formula and its Corollaries}\label{s1.2} In \cite{PV}
we stated the following theorem.

\begin{thm}[Theorem 1 of \cite{PV}]\label{T1} If $G$ is any based
Gauss diagram of a knot $K$, then
\begin{equation}\label{v2}
v_2(K)=\left\langle \vcenter{\epsffile{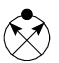}}, G\right\rangle.
\end{equation} \end{thm}

The right hand side is the sum $\sum\Ge(c_1)\Ge(c_2)$  over all
subdiagrams of $G$ isomorphic to $\vcenter{\epsffile{xup.eps}}$, where
 $c_1$, $c_2$ are the chords of the subdiagram.  General discussion on
the formulas of this kind see in \cite{PV}.

\begin{ex}\label{exv2} As it is easy to see from Figure
\ref{f2}, $v_2(4_1)=-1$.\end{ex}

\begin{cor}\label{Sym}If $G$ is any
based Gauss diagram of a knot $K$ then
\begin{equation}\label{v2sym}
v_2(K)=\left\langle \vcenter{\epsffile{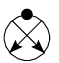}}, G\right\rangle.
\end{equation} \end{cor}

Corollary \ref{Sym} immediately follows from the fact that the 
rotation of the knot by $\pi$ around the $x$-axis results in a 
Gauss diagram of $K$ with all arrows of $G$ inverted, while their 
signs are preserved.

\begin{cor}\label{Arf}
If $G$ is a based Gauss diagram of a knot $K$, then the Arf 
invariant of $K$ is equal modulo 2 to the number of subdiagrams 
of $G$ isomorphic to $\vcenter{\epsffile{xup.eps}}$.
\end{cor}

There are a lot of methods for a calculation of the Casson knot 
invariant and the Arf invariant. See Kauffman \cite{Kauff}, Lannes 
\cite{Lannes}, \cite{Lannes'}, Gilmer \cite{Gilmer}.  
As far as we know, Theorems \ref{T1} and \ref{Arf} provide the 
easiest and the most practical ways for calculating $v_2$ and 
the Arf invariant.
The proof of Theorem \ref{T1} given below follows the lines of 
Kauffman's algorithm \cite{Kauff} for calculation of $\GD''_K(1)$.

\begin{thm}\label{Estimate} For any knot $K$ which admits
a diagram with $n$ crossing points $|v_2(K)|\le[\frac{n^2}8]$.
\end{thm}

\begin{rem}\label{rem1} Lin and Wang \cite{LW} found an estimate
of $|v_2(K)|$ which is twice weaker than Theorem \ref{Estimate}.
\end{rem}
\begin{rem}\label{rem2} The estimate of Theorem \ref{Estimate} is 
sharp for the case of odd $n$:  if $K$ is the torus knot of the type 
$(n,2)$ with any odd $n\ge3$ (it has $n$ crossings), then $v_2(K)$ 
is equal to $\frac{n^2-1}8$.  
For the case of even $n$ the inequality of Theorem \ref{Estimate} may 
be strengthened. 
An obvious consideration that the number of chords intersecting the 
given one is always even, permits to decrease the estimate at least 
by 1, but most probably this can be improved further. 
It is interesting whether the inequality $v_2\ge-[\frac{n^2}8]$, 
which follows from Theorem \ref{Estimate}, can be strengthened. 
\end{rem}

\subsection{The Proof of Theorem \ref{T1}}\label{s1.3} 
We use the following skein relation for the Casson knot invariant:
\begin{equation}\label{skeinforv2}
v_2(\vcenter{\epsffile{poscros.eps}})-v_2(\vcenter{\epsffile{negcros.eps}})
=lk(\vcenter{\epsffile{smtcros.eps}}).\end{equation}
In this formula, following a tradition, we present links (and knots)
by their fragments, which contain differences from other links under
consideration. By $lk(\vcenter{\epsffile{smtcros.eps}})$ it is denoted
the linking number of the components of link
$\vcenter{\epsffile{smtcros.eps}}$.

The relation \eqref{skeinforv2} is well-known. It is easy to check
that together with the condition $v_2(unknot)=0$ it defines a knot
invariant, this invariant is of degree 2 and takes the value 1 on 
trefoil.

To calculate $v_2$ of the knot $K$ presented by a diagram $G$, we
transform $K$ to the unknot, going from the base point along the 
orientation of $K$ and replacing an undercrossing by an 
overcrossing, if at the first passage through the point we go 
along the undercrossing. 
When we pass over the whole diagram, it becomes descending, and 
hence represents the unknot.
Each time we change a crossing $s$, the value of $v_2$ changes by
$-\Ge(s)lk(\vcenter{\epsffile{smtcros.eps}})$, where $\Ge(s)$ is the
sign of the crossing. Since $v_2(unknot)=0$, it gives
\begin{equation}\label{v2middle}
v_2(K)=\sum\Ge(s)lk(L_s),\end{equation}
where $L_s$ runs over links which appeared as smoothings at points 
where the crossing changed.

To calculate $lk(L_s)$, we can sum up the signs of all the crossing
points of $L_s$ in which the component containing the base point goes
below the other component. These points correspond to chords of $G$
intersecting the chord $c(s)$ corresponding to $s$ and directed to 
the side of $c(s)$ containing the base point. At the moment all arrows
of the original diagram $G$ with heads between the base point and
the head of $c(s)$ have been inverted. Therefore $lk(L_s)$ is equal 
to the sum of signs of arrows crossing $c(s)$ and having heads between
tail of $c(s)$ and the base point. In other words, $lk(L_s)$ is
$\sum\Ge(c_2)$ where the summation runs over all chords involved,
together with $c(s)$, into subdiagrams of the type
$\vcenter{\epsffile{xup.eps}}$.

Substituting this to \eqref{v2middle} we obtain \eqref{v2}.\qed

\subsection{Proof of Corollary \ref{Estimate}}\label{s1.4} Let $G$ be 
a based Gauss diagram of $K$ with $n$ chords. Subdivide the set $C$ of
all chords of $G$ into two subsets $C^+$ and $C^-$, where $C^+$ and
$C^-$ consist of all chords of the type $\vcenter{\epsffile{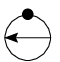}}$
and $\vcenter{\epsffile{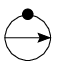}}$, respectively.  Let $|C^+|=k$,
$|C^-|=n-k$, $0\le k\le n$ and let $n_1$ and $n_2$ be the number of
subdiagrams of $G$ isomorphic to $\vcenter{\epsffile{xup.eps}}$ and
$\vcenter{\epsffile{xdown.eps}}$, respectively.  Any subdiagram of $G$
isomorphic to $\vcenter{\epsffile{xup.eps}}$ or
$\vcenter{\epsffile{xdown.eps}}$ consists of one chord from each 
subset $S^\pm$, thus $n_1+n_2\le k(n-k)$. It remains to notice that, 
as follows from Theorem \ref{T1} and Corollary \ref{Sym},
\begin{equation}v_2\le\min\{n_1,n_2\}\le\frac{k(n-k)}2\le
\left[\frac{n^2}8\right].
\end{equation}
\qed

\subsection{An Elementary Theory of the Casson Knot Invariant}
\label{s1.5}
Formula \eqref{v2} provides an elementary way to introduce the Casson
knot invariant. This formula, being used as a definition, gives a 
numeric function of a knot diagram with a marked point. 
At first glance, it is not clear if this is invariant with respect to
the isotopy.  However this is not difficult to check.

First, let us prove that
$\left\langle \vcenter{\epsffile{xup.eps}}, G\right\rangle$ does not
depend on the base point. When the base point moves along the circle 
of $G$, the expression
$\left\langle \vcenter{\epsffile{xup.eps}}, G\right\rangle$ can change
only at the moment of passing through an arrowhead. Denote this
arrow by $c$. Right before this moment the terms involving $c$ equal
to the product of $\Ge(c)$ by the sum of signs of all arrows 
crossing $c$ in the same direction.  Right after this moment, these 
terms are replaced by the product of $\Ge(c)$ by the sum of signs
of all arrows crossing $c$ in the opposite direction. Therefore to 
prove independence of the right hand side of \eqref{v2} on the base
point it suffices to notice, that for each chord $c$ of the Gauss 
diagram, the sum of signs of all arrows of the Gauss diagram crossing 
$c$ in one direction, is equal to the sum of signs of arrows crossing 
$c$ in the opposite direction. Indeed, these sums are equal to the 
linking number of the two-component link, obtained by smoothening the 
double point corresponding to $c$.

\begin{figure}[t]
\centerline{\epsffile{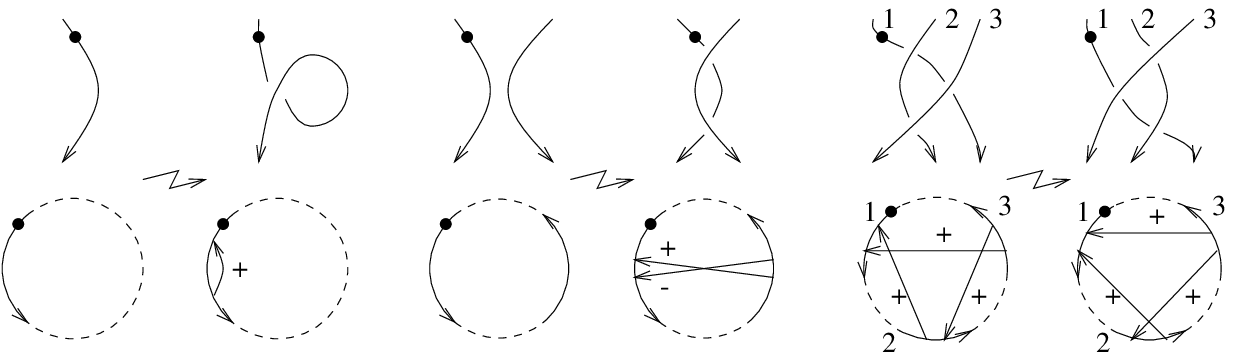}}
\caption{}
\label{Reidem}
\end{figure}

The invariance under Reidemeister moves follows from the study of
the corresponding changes of a Gauss diagram. See Figure 
\ref{Reidem}, where some of these changes are shown.

Under the first and third moves subdiagrams isomorphic to
$\vcenter{\epsffile{xup.eps}}$ do not change. Under the second move
each new subdiagram of this type includes exactly one of the two new
chords. Therefore new subdiagrams cancel out in pairs.

\section{From a Combinatorial Formula to Degrees of Maps}\label{s2}

\subsection{The Motivation: the Linking Number}\label{s2.00} Formula 
\eqref{v2} is similar to a combinatorial formula for linking number.  
Recall that the linking number of disjoint oriented circles 
$L_1,L_2\subset \R^3$ is equal to the sum of signs of all crossings 
in a diagram of the link $L_1\cup L_2$, where the $L_1$ passes over 
$L_2$. As we know, this can be interpreted as the formula for a 
calculation of the degree of map
\begin{equation}\label{Gf}\Gf:L_1\times L_2\to
S^2:(x,y)\mapsto\frac{x-y}{|x-y|}.\end{equation}
This suggests to look for a similar interpretation of \eqref{v2}.

\subsection{An Interpretation of the Casson Invariant via
a Local Degree}\label{s2.1}
Each summand in the expression for the right hand side of \eqref{v2} 
is a local degree for a map which is constructed as follows. 
For a knot $K\subset \R^3$ with a base point $*\in K$ denote by 
$C_{X}$ the space of 4-tuples $(x_1,x_2,x_3,x_4)\in K^4$ of points
ordered in the natural way, defined by the orientation of $K$:  
when one goes on $K$ along the orientation, the points occur in the 
sequence $*, x_1,x_2,x_3,x_4,*$.
Denote by $C^0_{X}$ the subspace of $C_{X}$ defined by inequalities
$*\ne x_1\ne x_2\ne x_3\ne x_4\ne *$. 
The orientation of $K$ and the order of coordinates determine an 
orientation of the manifold $C^0_{X}$. Define a map
$$\Gf^0_{X}:C^0_{X}\to S^2\times S^2:  (x_1,x_2,x_3,x_4)\mapsto
(\frac{x_1-x_3}{|x_1-x_3|},\frac{x_4-x_2}{|x_4-x_2|}).$$
The map $\Gf^0_{X}$ extends uniquely to the whole $C_{X}$ by 
continuity. Denote the extension by $\Gf_{X}$.

In the notations of the preceding paragraph $X$ stands for our
picture $\vcenter{\epsffile{xup.eps}}$.

The preimage of the point
$$(s,s)=(\text{south pole}, \text{ south pole})\in S^2\times S^2$$
under $\Gf_{X}$ consists of configurations of points corresponding to
subdiagrams  of the Gauss diagram isomorphic to
$\vcenter{\epsffile{xup.eps}}$. The contribution of a subdiagram to the
right hand side of \eqref{v2} is equal to the local degree of $\Gf_{X}$ at
the corresponding point of the preimage. Indeed, $\Gf_{X}$ is locally
equivalent to the Cartesian product of two copies of the map $\Gf$,
defined above by \eqref{Gf}. On the other hand, the local degree of
$\Gf$ is the sign of the corresponding chord.

Therefore $v_2(K)$ seems to be the degree of the map
$\Gf_{X}:C_{X}\to S^2\times S^2$. However, we have to be cautious:
the source space $C_{X}$ of this map is not a closed manifold. It is 
a manifold with boundary and corners.
Of course, we still can give a homology interpretation of the degree 
taking for the source space $C_{X}$ the relative homology 
$H_4(C_{X},C_{X}\sminus C^0_{X})$ and for the target $S^2\times S^2$ 
the relative homology 
$H_4(S^2\times S^2, S^2\times S^2\smallsetminus(s,s))$.

\subsection{We Run into Problems}\label{s2.2}This interpretation works
as long as the knot is in a position such that its vertical projection
is generic, i.e. gives rise to a knot diagram and the base point does
not coincide with a double point.  In particular, the projection is an
immersion without triple points and points of self-tangency.

However, the usual property of a degree to be invariant under
deformations does not follow and requires separate considerations.
Indeed, a generic knot isotopy involves situations when the projection
is not generic. These are exactly the moments when either the base
point passes the double point or the diagram experiences
Reidemeister moves. Let us treat these problems separately.

\subsection{Exiling Base Point to Infinity}\label{s2.3} The moments,
when the base point passes a double point, correspond to two
3-dimensional faces of our 4-dimensional space $C_{X}$. One of them
consists of configurations with $*=x_1\ne x_2\ne x_3\ne x_4\ne *$.
The other one is defined by $*\ne x_1\ne x_2\ne x_3\ne x_4= *$.
Denote them by $\GS X_{1*}$ and $\GS X_{4*}$, respectively.

Although there exists a natural homeomorphism
$(x_1,x_2,x_3,x_4)\mapsto(x_2,x_3,x_4,x_1)$ between them, we cannot
kill $\GS X_{1*}$ and $\GS X_{4*}$ by gluing via this homeomorphism, 
since it does not commute with the map $\Gf_{X}$.

To overcome this problem we resort to an old trick, which was used 
for example by Vassiliev \cite{V} in his original definition of 
Vassiliev knot invariants: to place the base point at the infinity. 
A $C^2$-smooth knot in $S^3$ with a base point is mapped by a
stereographic projection from the base point to a smooth knotted 
line in $\R^3$ with an asymptote. Moreover, by an arbitrary small 
diffeotopy one can turn a neighborhood of the base point on the 
original knot into a geodesic. This turns the image of the knot 
into a {\it long knot,\/} i.e., it coincides with a (straight) line 
outside of some ball. Without a loss of generality, we will assume 
this line to be the $y$-axis.

Since we need the space $C_{X}$ to be compact, in the case of long
knot $K$ it is constructed in a slightly different way.
First we compactify $K$ by adding a point at infinity. Denote the
compactified $K$ by $\tilde K$. Set $C_{X}$ to be the closure in
$\tilde K^4$ of $C^0_{X}\subset K^4$. Then extend $\Gf^0_{X}$ to
$C_{X}$ by continuity, as above.

This solves our problem: although the faces $\GS X_{1*}$ and
$\GS X_{4*}$ of $C_{X}$ (consisting of points with
$\infty=*=x_1\ne x_2\ne x_3\ne x_4\ne *$ and
$\infty=*\ne x_1\ne x_2\ne x_3\ne x_4= *$, respectively)
are still 3-dimensional, their images under $\Gf_{X}$ are
2-dimensional. Thus from homological point of view they are
unessential.

\subsection{Losing and Recovering the Degree During the Reidemeister
Moves}\label{s2.4} Consider now the strata of the boundary of $C_{X}$
which manifest themselves at Reidemeister moves.
For instance, at the third Reidemeister move (i.e. when the 
projection has a triple point), the isotopy is not a proper map 
over $(s,s)$ (which means that the preimage of $(s,s)$ meets the 
boundary). In other words some points of the preimage of $(s,s)$ jump
out of $C^0_{X}$ for an instant. The standard theory of degree based
on relative homology is designed for proper maps, and we cannot use 
it. One could hope that this happened because of a wrong choice of the
point and the situation can be improved by shifting $(s,s)$ off the
diagonal of $S^2\times S^2$.

\begin{figure}[htb]
$\centerline{\epsffile{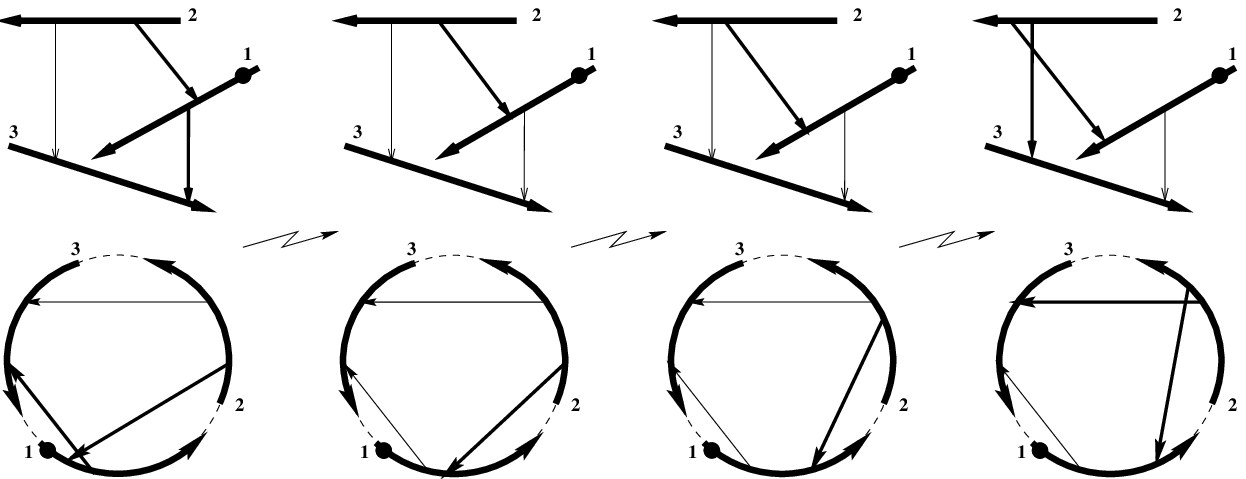}}$
\caption{}
\label{triple}
\end{figure}

However, instead of being improved the situation is getting even worse:
points of the preimage may not only appear on the boundary of
$C_{X}$ for an instant, but disappear for a certain period of time.
During this period the degree may jump several times.  In Figures
\ref{triple} we show how it happens. Three chords participate in 
this interaction. End points of two chords, involved in a subdiagram 
of the type $\vcenter{\epsffile{xup.eps}}$, meet and pass through each
other. The chords become disjoint. But then the opposite process
occurs, with another pair of chords. At that moment the original
degree of $\Gf_{X}$ is recovered. This suggests to look for a place 
where the degree was hidden.

\section{From a Local to a Global Degree}\label{s3}

\subsection{A Route from a Local to a Global Degree}\label{s2.6} 
\begin{figure}[htb]
$\centerline{\epsffile{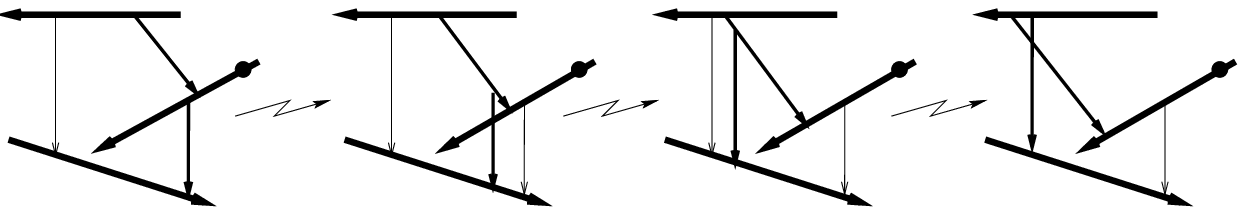}}$
\caption{}
\label{tripgraph}
\end{figure}

There is a nice solution of this puzzle: the chord which is involved 
in both pairs serves as a bridge between the point, where the first 
pair gets out of the game, and the point, where the second pair 
comes, and the second chord of the first pair may glide over this 
bridge. See Figure \ref{tripgraph}. 
On the way, there is a configuration of two oriented segments 
parallel to the fixed directions. One of the segments connects 
points on the knot, while the other one connects a point of the 
knot with a point on the first segment.

This resembles 3-valent graphs appearing in the Chern-Simons theory
approach to Vassiliev invariants, see Bar-Natan \cite{Bar-Nat} and 
Bott-Taubes \cite{BT}. Inspired by this picture, we combine our 
approach with the Chern-Simons approach below in this section.
More literally the same picture is used in Section \ref{s5.2}.

In the forthcoming sections we construct another configuration space
$C_{Y}$ related to a long knot and a continuous mapping 
$\Gf_{Y}:C_{Y}\to S^2\times S^2\times S^2$. Then we glue six copies 
of $C_{Y}$ and six copies of $C_{X}\times S^2$ together into a 
single 6-dimensional stratified space $\cal C$ with a fundamental 
class $[\cal C]\in H_6(\cal C, \cal S)$ where $\cal S$ is a union 
of some low-dimensional strata of $\cal C$. The maps 
$\Gf_{X}\times \id_{S^2}$ and $\Gf_{Y}$ give rise to a continuous
map $\Gf:\cal C\to S^2\times S^2\times S^2$, which maps $\cal S$ 
into a 5-dimensional set $D\subset S^2\times S^2\times S^2$, 
consisting of triples $(u_1,u_2,u_3)$ of coplanar vectors which are 
either coplanar or contain the vector $(0,\pm1,0)\in S^2\subset\R^3$ 
(the latter is due to our convention that long knots coincide with 
the $y$-axis at the infinity).
We prove that the degree of this map is
$6v_2(K)$ (i.e., $\Gf[\cal C]=6v_2(K)[S^2\times S^2\times S^2]$, 
where $[\cal C]$ is the natural generator of $H_6(\cal C,\cal S)$ 
and $[S^2\times S^2\times S^2]$ is the orientation generator of
$H_6(S^2\times S^2\times S^2)$).

\subsection{The Principal Faces of $C_{X}$}\label{s2.5} Prior to
construction of a space completing $C_{X}$ to a cycle, we have to
enlist the codimension one faces of $C_{X}$. Two of them, $\GS X_{1*}$
and $\GS X_{4*}$, were studied above in Section \ref{s2.3}. Passing
to long knots made them unessential.

The other principal strata of the boundary can be described as follows.
For any $i\in\{1,2,3\}$ denote by $\GS X_{i\ i+1}$ the subset of
$C_{X}$ consisting of points $\Gg\in C_{X}$, such that the components 
$x_{j}$  of $\pi(\Gg)$ are distinct from $*=\infty$, and distinct 
from each other except for $x_i=x_{i+1}$.

\begin{lem}\label{2.1.principle_faces} Let $A$ be one of the pairs 
$1*$, $12$,  $23$, $34$, or $4*$.  Then the space $\GS X_{A}$ is a 
3-dimensional manifold open in $C_{X}\smallsetminus C^0_{X}$.
The union $C^0_{X}\cup\GS X_{A}$ is a manifold with boundary
$\GS X_{A}$.  The complement of $\bigcup_{A}\GS X_A$ in
$C_{X}\smallsetminus C^0_{X}$ has dimension 2.\end{lem}

This is a straightforward consequence of the construction of $C_{X}$.
In fact, one can identify $C_{X}$ with a closed 4-simplex. In the
coordinates $x_1$, $x_2$, $x_3$, $x_4$ in $C_{X}$ the strata of 
$C_{X}\smallsetminus C^0_{X}$ are described by obvious linear 
equations and inequalities.
Note that the faces $\GS X_{12}$, $\GS X_{23}$, $\GS X_{34}$ are
naturally homeomorphic to the same space. Namely, denote by 
$C^0_{V}$ the space of 3-tuples $(x_1,x_2,x_3)\in K^3$ of points 
ordered in the natural way defined by the orientation of $K$ with 
$*\ne x_1\ne x_2\ne x_3\ne *$. This is a 3-dimensional manifold
equipped with the orientation determined by the order of coordinates
and the orientation of $K$.

\begin{lem}\label{2.2.xi} The maps
$$\begin{aligned}&\xi_1:\GS X_{12}\to
C^0_{V}:(x_1,x_2,x_3,x_4)\mapsto(x_2,x_3,x_4)\\
&\xi_2:\GS X_{23}\to C^0_{V}:(x_1,x_2,x_3,x_4)\mapsto(x_1,x_3,x_4)\\
&\xi_3:\GS X_{34}\to
C^0_{V}:(x_1,x_2,x_3,x_4)\mapsto(x_1,x_2,x_3)\end{aligned}$$
are homeomorphisms.
The degree of $\xi_i$, $i=1,2,3$ with respect to the orientation
induced on $\GS_{i\ i+1}$ as on the boundary of $\GG\cup\GS_{0i}$
is $(-1)^i$.
\end{lem}

\subsection{The Configuration Space $C_{Y}$}\label{s2.7} Consider the
space $C^0_{Y}$ of 4-tuples
$$(x_1,x_2,x_3,x_0)\in K^3\times\R^3,$$
where $x_1$, $x_2$, $x_3$   are distinct from each other and $x_0$
and ordered in the natural way which is determined by the orientation
of $K$. Here in the notations $Y$ represents picture
$\vcenter{\epsffile{ygraph.eps}}$.
The order of the coordinates, and the orientations of $K$ and
$\R^3$, determine an orientation of the manifold $C^0_{Y}$.
Define a map $$\Gf^0_{Y}:C^0_{Y}\to S^2\times S^2\times S^2:
(x_1,x_2,x_3,x_0)\mapsto\left(\frac{x_1-x_0}{|x_1-x_0|},
\frac{x_3-x_0}{|x_3-x_0|},\frac{x_0-x_2}{|x_0-x_2|}\right).$$

As in the previous section, we need to embed $C^0_{Y}$ to some
compact space and extend $\Gf^0_{Y}$ to it. The former is easy:
$C^0_{Y}\subset \tilde K^3\times S^3$, but $\Gf^0_{Y}$ does not
admit a continuous extension to $\tilde K^3\times S^3$. We may use a
standard way of overcoming this difficulty: to consider the graph $\GG$
of $\Gf^0_{Y}$ as a subset of
$\tilde K^3\times S^3\times S^2\times S^2\times S^2$
and take the closure.  Denote the closure by $C_{Y}$
and its image under the natural projection
$\pi:C_{Y}\to\tilde K^3\times S^3$ by $B_{Y}$. This is a sort of a
resolution of singularities: the restriction to $\GG$ of
$\pi:C_{Y}\to B_{Y}$ identifies $\GG$ with $C^0_{Y}$.  Via this
identification the natural projection $C_{Y}\to S^2\times S^2\times
S^2$ extends the original map $\Gf^0_{Y}$. Denote this extension by
$\Gf_{Y}$.

Our space $C_{Y}$ can be identified with a subspace of a quotient
space of the widely-known and well-studied space $C_{3,1}$, obtained 
from $C^0_{Y}$ by the Fulton-MacPherson \cite{FM} compactification 
construction. (Similarly, $C_{X}$ is a quotient space of the space 
$C_4$ obtained by an analogous compactification of $C^0_{X}$.) 
Various aspects of this construction were presented with details in 
\cite{AS} and \cite{BT}.  The difference  between $C_{Y}$ and a space 
studied in \cite{BT} is that we study a long knot (or a knot in $S^3$ 
with a base point), and, furthermore, we make the minimal resolution 
of singularities needed to define $\Gf_{Y}$, while Bott and Taubes 
use a larger Fulton-MacPherson compactification \cite{FM} of the 
configuration space.  
Thus our space $C_{Y}$ turns out to be a subspace of a quotient space 
of $C_{3,1}$ from \cite{BT}.  However, for our purposes we do not 
need a refined analysis of the natural stratification of $C_{3,1}$ 
presented in \cite{BT}.  Instead, we use the following elementary 
consideration of the boundary $C_{Y}\smallsetminus\GG$ of $\GG$.

\subsection{Principal Faces of $C_{Y}$}\label{s2.8}
Let $\Gg=(x_1,x_2,x_3,x_0,u_1,u_2,u_3)$ be a point of
$C_{Y}\smallsetminus\GG$.  Since
$$\pi(\Gg)=(x_1,x_2,x_3,x_0)\in B_{Y}$$ belongs to the boundary of
$C^0_{Y}$, either $x_0=\infty$ $(= S^3\smallsetminus\R^3)$, or 
$x_1=*=\infty$, or $x_3=*=\infty$, or some of $x_i$ coincide. 
Consider all the cases separately.

For any subset $A$ of $\{0,1,2,3,*\}$ denote by $\GS Y_{A}$ the subset
of $C_{Y}$ consisting of points
$\Gg=(x_0, x_1,x_2,x_3,u_2,u_2,u_3)\in C_{Y}$ such that in the
configuration  $x_0, x_1,x_2,x_3,*$  two points coincide if
and only if the corresponding elements of $\{0,1,2,3,*\}$ belong to $A$.
For instance,
$$\GS Y_{01*}=\{\Gg\in C_{Y}\,:\,*= x_0=x_1\ne x_2\ne x_3\ne *\}.$$
Of course, $\GS Y_{A}$ with $A=\{1,3\}$ and $\{2, *\}$ are empty.

Several boundary strata $\GS Y_A$ are of codimension 2 or higher in
$C_{Y}\smallsetminus\GG$.

\begin{lem}\label{2.8.small_strata}
Let $A$ be a subset of $\{1,2,3,*\}$ containing at least 3 elements.
Then $\dim(\GS Y_A)\le 4$.\end{lem}

\begin{proof}
Since $0\notin A$, the map $\Gf^0_{Y}$ extends uniquely to
$\Gp(\GS Y_A)\subset B_{Y}$ by  continuity.
Thus, by the construction of $C_{Y}$, the stratum $\GS Y_A$ is projected
homeomorphically to $\Gp(\GS Y_A)$.
As the codimension of $\Gp(\GS Y_A)$ in $B_{Y}$ is $|A|-1\ge 2$,
$\dim(\GS Y_A)=\dim(\Gp(\GS Y_A))\le 4$.
\end{proof}

 The rest of non-empty boundary strata are of codimension 1 in
$C_{Y}\smallsetminus\GG$.
The strata $\GS Y_{0i}$, $i=1,2,3$ are of primary interest,
as, similarly to $\GS X_{i\ i+1}$, they are homeomorphic to
$C^0_{V}\times S^2$.

\begin{lem}\label{2.8.eta} For $i=1,2,3$, the map
$$\eta_i:\GS Y_{0i}\to C^0_{V}\times S^2:
(x_1,x_2,x_3,x_0,u_1,u_2,u_3)\mapsto((x_1,x_2,x_3),u_i)$$ is a
homeomorphism of degree $(-1)^i$ with respect to the orientation
induced on $\GS Y_{0i}$ as on the boundary of
$\GG\cup\GS Y_{0i}$ and the product of the orientation of
$C^0_{V}$ defined above by the standard orientation of $S^2$.
\end{lem}

Some other strata, which seem to be rather big, admit orientation
reversing homeomorphisms. In the next section this allows us to
cancel them out.

\begin{lem}\label{2.2.E}
Let $A=\{0,1,2\}$ or $\{0,2,3\}$. The stratum $\GS Y_A$ is
a codimension 1 submanifold of a manifold $C^0_{Y}\cup\GS Y_A$.
The maps
$$\Gz_1:\GS Y_{012}\to\GS
Y_{012}:
(x_1,x_2,x_3,x_0,u_1,u_2,u_3)\mapsto
(x_1,x_2,x_3,x_0,u_2,u_1,u_3), $$
$$\Gz_2:\GS Y_{023}\to\GS
Y_{023}:
(x_1,x_2,x_3,x_0,u_1,u_2,u_3)\mapsto
(x_1,x_2,x_3,x_0,u_1,u_3,u_2) $$
are homeomorphisms which can be extended to
orientation reversing homeomorphisms of a neighborhood of
$\GS Y_A$ in $C^0_{Y}$.
\end{lem}
\begin{proof} The extensions can be defined by the following formulas:
$$
(x_1,x_2,x_3,x_0)\mapsto
(x_1,x_2,x_3,x_1+x_2-x_0), $$
$$
(x_1,x_2,x_3,x_0)\mapsto
(x_1,x_2,x_3,x_2+x_3-x_0) $$
\end{proof}
%

\subsection{Gluing Pieces Together}\label{s.2.9}
\begin{figure}
\centerline{\epsffile{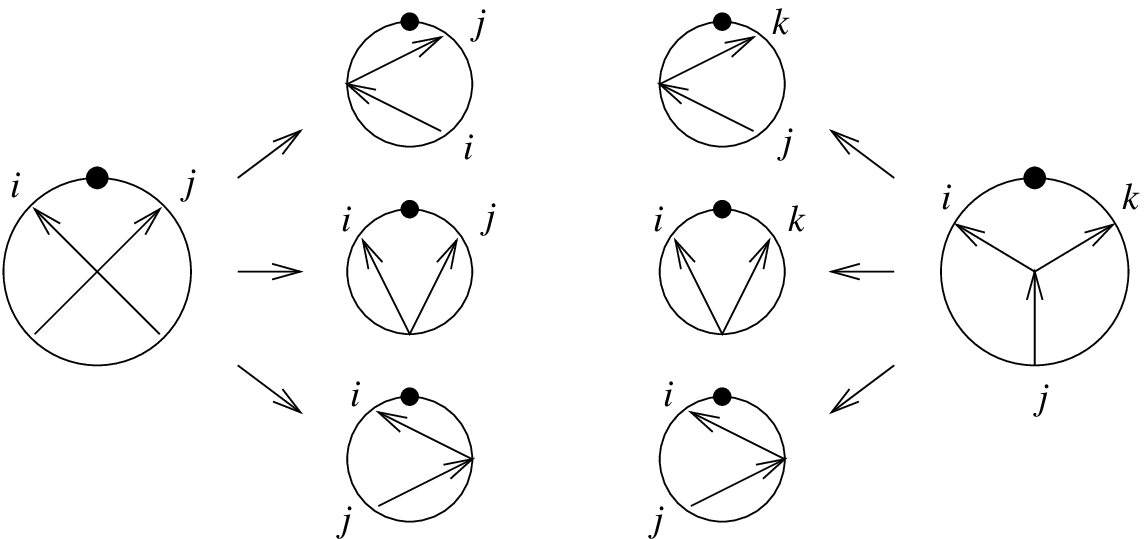}}
\caption{}
\label{gluing}
\end{figure}

Now we are to construct $\cal C$ as outlined in
Section \ref{s2.6}. We consider 6 copies of $C_{X}\times S^2$ and
$C_{Y}$, i.e. the product $(C_{X}\times S^2\cup C_{Y})\times S_3$.
Here the symmetric group $S_3$ is equipped with the discrete topology.
The space $\cal C$ is obtained as the quotient space of
$(C_{X}\times S^2\cup C_{Y})\times S_3$ by the following
identifications. \begin{enumerate}
\item $\GS X_{12}\times S^2\times\Go$ is identified with
$\GS Y_{01}\times\Go\circ(1,3,2)$ via
$(\xi_1\times\id_{S^2})\circ\eta_1^{-1}$;
\item $\GS X_{23}\times S^2\times\Go$ is identified with
$\GS Y_{02}\times\Go\circ(2,3)$ via
$(\xi_2\times\id_{S^2})\circ\eta_2^{-1}$;
\item $\GS X_{34}\times S^2\times\Go$ is identified with
$\GS Y_{01}\times\Go$ via
$(\xi_3\times\id_{S^2})\circ\eta_3^{-1}$;
\item $\GS Y_{012}\times\Go $  is identified with
$\GS Y_{012}\times \Go\circ(1,2)$ via $\Gz_1$;
\item $\GS Y_{023}\times\Go $  is identified with
$\GS Y_{023}\times \Go\circ(2,3)$ via $\Gz_2$;
\item the induced identifications on the boundaries of the strata above
\end{enumerate}
For $\Go\in S_3$, let $\bar\Go:S^2\times S^2\times S^2\to
S^2\times S^2\times S^2$ be the permutation of the
factors defined by $\Go$.
One can easily check that, as it was promised above, the maps
$\bar\Go\circ\Gf_{X}\times \id_{S^2}:C_X\times S^2\times\Go\to
S^2\times S^2\times S^2$ and $\bar\Go\circ\Gf_{Y}:C_Y\times\Go\to
S^2\times S^2\times S^2$ give rise to a continuous
map $\Gf:\cal C\to S^2\times S^2\times S^2$. See Figure \ref{gluing}.

Despite of all these identifications, $\cal C$ is not closed yet,
that is although its high-dimensional strata are orientable and gluing
reverses orientations, the high-dimensional homology group $H_6(C)$ is
trivial. The boundary five-dimensional strata are obtained from
$\GS Y_A$ with $A\subset\{0,1,2,3,*\}$ containing both $0$ and $*$ or
a pair of consecutive elements of the sequence $\{*,1,2,3,*\}$, or
$A=\{0,1,2,3\}$.
Denote the union of these strata by $\cal S$. It is easy to see that
$H_6(\cal C,\cal S)=\Z$. Indeed, the six-dimensional strata
$C^0_{Y}\times \Go$ and $C^0_{X}\times S^2\times\Go$ are connected and
oriented. They are attached to each other by orientation reversing
diffeomorphisms of five-dimensional strata on their boundary into a
connected space $\cal C$. Finally, $\cal S$ is the union of all the
five-dimensional strata not involved in the gluing.

Now let us study the image of $\cal S$ under $\Gf$ and show that it is
contained in the set $D$ of triples $(u_1,u_2,u_3)\in S^2\times S^2\times S^2$
of vectors which are either coplanar or contain vector $\pm a$,  where
$a=(0,1,0)\in S^2\subset\R^3$. 

\begin{lem}\label{2.8.small_im} Let $A$ be a subset of 
$\{1,2,3,*\}$.  Then $\Gf_{Y}$ maps $\GS Y_A$ into $D$.  \end{lem} 

\begin{proof} Observe, that if $A$	contains two elements of 
$\{1,2,3,*\}$, then it contains a pair of consecutive elements of the 
sequence $\{*,1,2,3,*\}$. If $A\supset\{1,*\}$ or $A\supset\{3,*\}$, 
then $u_1=-a$ or $u_3=a$, respectively.
If $A\supset\{1,2\}$ or $A\supset\{2,3\}$, then $u_1=-u_2$ or 
$u_2=-u_3$, respectively. \end{proof}

\begin{lem}\label{2.8.coplanar}
Let $A$ be a subset of $\{0,1,2,3,*\}$ containing both
$0$ and $*$.
Then $\Gf_{Y}$ maps $(\GS Y_A)$ into $D$.\end{lem}
\begin{proof} When $x_0$ tends to infinity, all three vectors $u_1$,
$u_2$, $u_3$ lie in the plane containing $a$ and the direction of the
move of $x_0$. 
\end{proof}

\begin{lem}\label{2.9}
Let $A=\{0,1,2,3\}$.
Then $\Gf_{Y}$ maps $(\GS Y_A)$ into $D$.\end{lem}
\begin{proof} All three vectors $u_1$, $u_2$, $u_3$ lie in the plane containing the direction
of the tangent vector to $K$ at $x_0=x_1=x_2=x_3$.
\end{proof}

\begin{thm}\label{main_th} The space $\cal C$	has a
well-defined fundamental class $[\cal C]\in H_6(\cal C,\cal S)$.
The map
$\Gf:\cal C\to S^2\times S^2\times S^2$	induces homomorphism
$H_6(\cal C,\cal S)\to H_6(S^2\times S^2\times S^2,D)$,
which maps $[\cal C]$ to $6v_2(K)[S^2\times S^2\times S^2]$. 
\end{thm}

\begin{proof}Lemmae \ref{2.8.small_im}- \ref{2.9} prove the first
statement of Theorem \ref{main_th}. Now we prove the rest.

To evaluate the degree, we return to the arguments given in the first
subsection of this section. Assume that our knot $K$ is in general
position with respect to the vertical projection. Calculate the degree
by counting (with signs) points of the preimage of a regular value
$r\not\in D$ of $\Gf:\cal C\to S^2\times S^2\times S^2$ close  to
$(s,s,s)$.
 As we observed in Section \ref{s2.1}, those of them which
belong to each of the six copies of $C_{X_*}\times S^2$ contribute
$v_2(K)$. It remains to notice that the preimage does not intersect
the copies of $C_{Y_*}$.
Indeed, each point of the preimage belonging to one of these copies
would correspond to a configuration of $(x_1, x_2, x_3)\in K^3$  such
that $x_2$ is positioned in $\R^3$ almost strictly above $x_1$ and
$x_3$. The points $x_1$ and $x_3$ are not close to each other because
they are separated on $K$ by $x_2$. The projection of $K$ is assumed to
be generic. In particular, it does not have triple points. Therefore
if $r$ is sufficiently close to $(s,s,s)$, this configuration cannot
appear.
Note also, that although $D$ divides $S^2\times S^2\times S^2$, and the
regular value $r$ may be chosen in any component of the complement of
$D$, the above evaluation of the local degree does not depend on this
choice.
This completes the proof of Theorem \ref{main_th}.\end{proof}

\subsection{Straightforward Applications}\label{s.straight}
There are different ways for calculating the degree of a
map. Two of them are classical.

First, one can take a regular value of the map and count the
points of its preimage with signs (which are the local degrees of this
map at the point). For instance, choosing a  point sufficiently close to
$(s,s,s)\in S^2\times S^2\times S^2$ we get again our combinatorial
formula \eqref{v2}. Choosing a point sufficiently close to 
$(-s,-s,-s)\in S^2\times S^2\times S^2$ we get \eqref{v2sym}.
Other choices of regular values give rise to other completely different
combinatorial formulas for $v_2$. Indeed, the strata of $\cal C$, which
are copies of $C_Y$, have been added just to make the space closed, but
under the original choice of the regular value they did not give any
input in the combinatorial formulas. However under other choices of the
regular value they become visible and change the type of the formulas.
We will deal with this in the next section.

The second classical way is to take a differential form of
top degree on the target, normalize it by the condition that the
integral of this form over the whole manifold is equal to one, pull it
back to the source and integrate over the whole source manifold.
In this way one can deduce from Theorem \ref{main_th} Bar-Natan's
integral formula \cite{Bar-Nat} for $v_2(K)$.

These two methods can be mixed which gives rise to a method 
generalizing both of them. See Section \ref{s5.1} below.

\section{From Regular Values to New Combinatorial Formulae}\label{s4}

\subsection{Counting Tinkertoy Diagrams}\label{s4.1}

If we choose a generic regular value of $\Gf$ and do not impose any
restriction on the position of a knot, counting preimages of the
the regular value reduces to counting configurations of arrows in the
space of the following two types. The configurations of the first type
are pairs of arrows connecting points of the knot. The arrows are
attached to the knot according to arrow diagram
$\vcenter{\epsffile{xup.eps}}$ and are directed in two of the three
fixed directions. The configurations of the second type are tripods
made of three arrows connecting a point in the space with three points
on the knot. The arrows are attached according to the diagram
$\vcenter{\epsffile{ygraph.eps}}$ and are directed in the three fixed
directions.

Similar configurations have been considered by D. Thurston
\cite{Thurston} under the name of {\it tinkertoy diagrams.\/} The only
difference is that his diagrams consist of unoriented segments, while
ours are made of arrows. So, we will use the same term {\it tinkertoy
diagram.\/}

Theorem \ref{main_th} implies curious geometric consequences concerning
numbers of various tinkertoy diagrams on a given knot. However, we do
not elaborate this topic in a full generality. Instead, we obtain new
combinatorial formulas related to special kinds of knot diagrams.

\subsection{Regular Values near Both Poles}\label{s4.1.1}
One of the interesting choices of the regular value is to take a point
close to $(s,s,-s)$. Recall that since the vectors $s$, $s$, $-s$ are
coplanar and hence $(s,s,-s)\in D$ , the point  $(s,s,-s)$ cannot be
used as a regular value of $\Gf$ to calculate $v_2(K)$ via Theorem
\ref{main_th}. The same happens with  $(s,s,s)$ and to get our
combinatorial formula \eqref{v2} we took a regular value of $\Gf$
sufficiently close to  $(s,s,s)$.  All points sufficiently
close to  $(s,s,s)$ give the same combinatorial formula.

However points close to  $(s,s,-s)$ give rise to different
combinatorial formulas. Since we are interested in limit situations, it
is reasonable to consider smooth paths $t\mapsto(s_1(t),s_2(t),s_3(t))$
with $\lim_{t\to\infty}(s_1(t),s_2(t),s_3(t))=(s,s,-s)$, check if the
numbers of tinkertoy diagrams stabilize after some value of $t$ and
write down the combinatorial formula obtained for sufficiently large
$t$. However, first, let us consider the tinkertoy diagrams
corresponding to a generic  point close to  $(s,s,-s)$.

\subsection{Pairs of Arrows}\label{s4.1.2}
Tinkertoy diagrams of the first kind (i.e., pairs of arrows with the
end-points on the knot) consist of almost vertical, i.e. almost
parallel to $z$-axis, arrows. Hence these arrows appear near double
points of the knot projection to $xy$-plane. However, not all of them
are directed downwards: one can be directed upwards. Therefore the
tinkertoy diagrams of the first kind appearing at pairs of double
points of the knot projection to $xy$-plane make a contribution
different from the contribution in the case of a regular value close to
$(s,s,s)$. Recall that then the contribution in the case of $(s,s,s)$
was just $\left\langle6\ \vcenter{\epsffile{xup.eps}}\ , G\right\rangle$.
Now it is 
$$\left\langle2\ \vcenter{\epsffile{xup.eps}}+2\ \vcenter{\epsffile{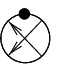
}}+2\ \vcenter{\epsffile{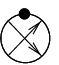}}\ , G\right\rangle.$$

\begin{figure}
\centerline{\epsffile{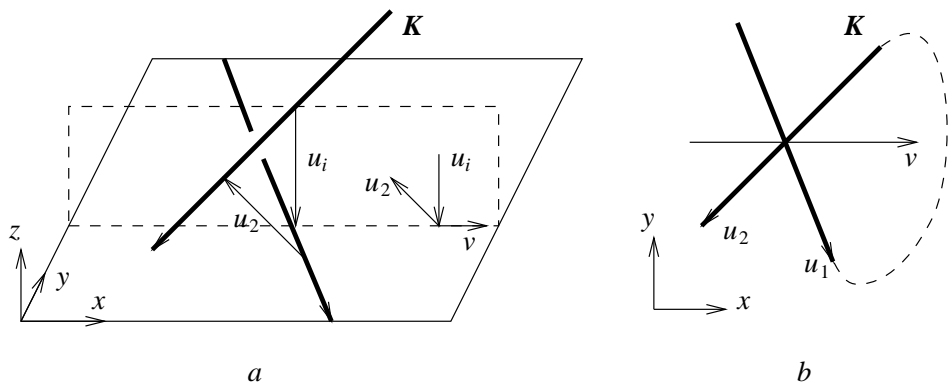}}
\caption{Tinkertoy diagram near a double point}
\label{contrib}
\end{figure}

A pair of almost vertical chords can be found also near the same double
point of the knot projection. On a tinkertoy diagram of the first type
the arrowheads are separated on the knot by the arrowtails (and by the
base point $*=\infty$ on the other side). So the arrowheads cannot be
close to each other on the knot. Therefore although a pair of almost
vertical chords can be found near the same double point of the knot
projection, in the case of $(s,s,s)$ tinkertoy diagrams of the first
type do not appear near the same double point. In the case of
$(s,s,-s)$ this may happen, see Figure \ref{contrib}a. A double point
$c$ where a tinkertoy diagram corresponding to a pair of vectors
$(s_i,s_3)$ with $i=1$ or $2$ appears, may be described by the
following combinatorial rule. Consider the plane $P$ spanned by $s_i$
and $s_3$ and passing through $c$. Choose a vector $v$ in the
intersection of $P$ with the plane of projection so that the
orientations of $P$ defined by the frames $(s_i,s_3)$ and $(s_3,t)$
coincide. Denote by $t_1$, $t_2$ the tangent vectors to the branches of
the knot projection in $c$ (oriented and ordered by the orientation of
our long knot). Then the condition is that the orientations of the
$xy$-plane of projection, defined by the frames $(v,t_1)$, $(v,t_2)$
and $(t_1,t_2)$ coincide (due to generic choice of $s_i$, $s_3$,
vectors $t_1$, $t_2$ are transversal to $P$). See Figure
\ref{contrib}b.

Thus to keep the contribution from the tinkertoy diagrams of the first
type fixed, one needs to keep planes spanned by $s_1$, $s_3$ and $s_2$,
$s_3$ unchanged as $(s_1,s_2,s_3)$ approaches $(s,s,-s)$.

\subsection{Triangles Inscribed in Knot Diagram}\label{s4.1.3}
 Tripod tinkertoy diagrams behave in a more complicated way. If the
regular value $(s_1,s_2,s_3)$ is chosen near $(s,s,-s)$ generically,
then the free vertex $x_0$ of the tripod is positioned high over the
knot, as a vertex of sharp triangular pyramid with corners on the knot.
The corner corresponding to $s_3$ should be between two other corners
as they appear along the knot. This corner will be referred as
{\it northern.\/}
If the knot was positioned in the plane of the projection, the pyramids
of this sort would correspond to triangles homothetic to each other and
inscribed in the knot diagram in such a way that on the knot the
northern vertex lies between two other vertices. If $(s_1,s_2,s_3)$ is
sufficiently close to $(s,s,-s)$ then the differences between the
tripod tinkertoy diagrams based on the knot and on the knot diagram
become inessential. In this case the contribution coming from tripods
depends only on the $xy$-plane projection of the knot.

Now we can fix a curve in $S^2\times S^2\times S^2$ approaching
$(s,s,-s)$ in such a way that the combinatorial formula for $v_2(K)$
defined via Theorem \ref{main_th} by counting tinkertoy diagrams
associated with the point of this curve stabilizes as the point
approaches $(s,s-s)$:

Choose a triangle $T$ in $xy$-plane generic with respect to the knot
diagram under consideration. The genericity here means that the sides
of $T$ are not parallel to the tangent lines to the branches of the
knot projection at double points. Let $(X_1,Y_1)$, $(X_2,Y_2)$,
$(X_3,Y_3)$ be vertices of $T$. Denote by $s_i(t)$ with $i=1,2,3$, $t>0$
the unit vectors $\frac{(X_i,Y_i,-t)}{\sqrt{X_i^2+Y_i^2+t^2}}$ for
$i=1,2$ and $-\frac{(X_3,Y_3,-t)}{\sqrt{X_3^2+Y_3^2+t^2}}$ for $i=3$.
This is a smooth curve with
$\lim_{t\to\infty}(s_1(t),s_2(t),s_3(t))=(s,s,-s)$. As we saw, the set
of tinkertoy diagrams associated with $(s_1(t),s_2(t),s_3(t))$
stabilizes as $t\to\infty$ and the resulting combinatorial formula for
$6v_2(K)$ depends only on the diagram of $K$. It contains three terms:
\begin{enumerate}
\item
$\left\langle2\ \vcenter{\epsffile{xup.eps}}+
2\ \vcenter{\epsffile{xleft.eps}}+
2\ \vcenter{\epsffile{xright.eps}}\ , G\right\rangle,$
\item the number of double points of the projection positioned
in the way described above with respect to the sides connecting the
northern vertex $(X_3,Y_3)$ with the other two vertices,
\item and the
algebraic number of triangles homothetic to $T$ and inscribed in the
knot projection in such a way that on the knot the northern vertex lies
between the two others.
\end{enumerate}
The triangles are counted with signs. An
interested reader can find a combinatorial rule for the sign of an
inscribed triangle. Of course, it is nothing but the local degree of
$\Gf$ at the corresponding point.

Triangles inscribed in a knot projection are not customary for knot
theory. Choosing more sophisticated paths approaching $(s,s,-s)$, we
will derive new combinatorial formulas involving more common
characteristics of knot projection.

\subsection{Degeneration of Triangles}\label{s4.1.4}
Fix a positive $\GD$ and consider a family of triangles $T_t$ in
$xy$-plane  with vertices $(-1,0)$, $(\GD,0)$, $(0,\frac1t)$. These are
triangles with the same base $[-1,\GD]$ and height tending to $0$ as
$t\to\infty$. Replace in the construction of the path
$t\to(s_1(t),s_2(t),s_3(t))$ the triangle $T$ with $T_t$: put
$s_1(t)=\frac{(-1,0,-t)}{\sqrt{1+t^2}}$,
$s_2(t)=\frac{(\GD,0,-t)}{\sqrt{\GD^2+t^2}}$ and
$s_3(t)=-\frac{(0,1/t,-t)}{\sqrt{t^2+1/t^2}}$.

Assume that at double points of the knot projection there is no branch
with  tangent parallel to $x$-axis. Then the tinkertoy diagrams
associated to the points of the path under consideration stabilizes as
$t\to\infty$. The diagrams of the first type look as in the previous
case. Because of the special choice of the triangles, the combinatorial
rule for calculating the number of double points with tinkertoy diagram
of the first type simplifies and gives the number of double points of
the knot projection  where both branches are oriented upwards or both
downwards.

The tripod tinkertoy diagrams are of two sorts. The ones of the first
sort are related to the points of the knot projection where the
$y$-coordinate restricted to the knot projection has a local maximum.
The corresponding inscribed triangle shrinks to this point. The
contribution to the formula is $-1$.

The tinkertoy diagrams of the second sort are related to triples of
points of the knot projection satisfying the following conditions.
The points lie on the same line parallel to $x$-axis. The ratio of
the distances between the middle point
and two end points equals $\GD$. The middle point arises from the
northern vertex, and hence lies between the
other two points both on this horizontal line and on the knot.
The contribution of such a triple is
$\Ge=\pm1$ defined by the following rule.

Denote the points of the triple by $a$, $b$, $c$ in order of their
appearance on the knot. Moving the horizontal line containing $(a,b,c)$
up, we include $(a,b,c)$ in a one-parameter family of triples
$(a(y),b(y),c(y))$ of points of the knot projection. Denote by $\Gs$
the sign of the derivative of $|c(y)-b(y)|-|a(y)-b(y)|$ at the initial
position. For example, on the left hand side of Figure \ref{GD-tangle}
$\Gs=-1$, while on the right hand side $\Gs=1$. Denote by $q$ the
number of branches of the knot projection passing through $a$, $b$, $c$
upwards. Then $\Ge=(-1)^q\Gs$.

The formula which is obtained in this way still involves geometry of
the knot, though the geometry is reduced to planar geometry of the knot
projection. Choosing an appropriate $\GD$ or deforming a knot diagram,
one can make the formulas purely combinatorial. We will do this in 
Sections \ref{s4.2}-\ref{s4.4}. 

\subsection{Other Paths to Degeneration}\label{s4.1.5}
Choosing other paths to $(s,s,-s)$ one can get many other combinatorial
formulas. Even a simple renumeration of the vertices of $T_t$ changes
the result. We consider two renumeration.

For the first of them, put 
$s_1(t)=\frac{(-1,0,-t)}{\sqrt{1+t^2}}$,
$s_2(t)=\frac{(0,1/t,-t)}{\sqrt{t^2+1/t^2}}$ and
$s_3(t)=-\frac{(\GD,0,-t)}{\sqrt{\GD^2+t^2}}$.
Literally repeating the arguments of Section \ref{s4.1.4} we have to
make the following changes.

First, the same combinatorial rule for calculating the number of double
points with tinkertoy diagram of the first type gives twice
the number of double points of $f$ where either both branches are
oriented upwards and their intersection number\footnote{To define the
intersection number, one needs orientation and order of branches. Both
are defined by the orientation of the source line $\R$ of the
immersion.} is $+1$, or  both branches are oriented downwards and the
intersection number is $-1$.

Second, counting the contribution from tripods we observe that the one 
of local maxima disappears. The reason is that we have to count only
inscribed triangles whose northern vertex lies on the knot between two
other vertices. In this case, the rightmost vertex is northern, so such
a triangle cannot be inscribed at the maximum respecting the order.

Third, by the same reason, triples of points of the knot projections
lying on the same horizontal line should appear in another order on the
knot: the rightmost point on the line should be the middle one on the
knot.

Another renumeration is provided by
$s_1(t)=\frac{(\GD,0,-t)}{\sqrt{\GD^2+t^2}}$,
$s_2(t)=\frac{(0,1/t,-t)}{\sqrt{t^2+1/t^2}}$ and
$s_3(t)=-\frac{(-1,0,-t)}{\sqrt{1+t^2}}$.
Similarly to the above, the contribution of double points gives twice 
the number of double points of $f$ where either both branches are 
oriented upwards and their intersection number is $-1$, or both 
branches are oriented downwards and the intersection number is $+1$.  
The contribution made by tripods comes from triples of points of the 
knot projections lying on the same horizontal line such that the 
leftmost point on the line is the middle one on the knot.  

\subsection{Regular and Nonassociative Immersions}\label{s4.2}
Now we have to make preparations for reformulating results in a purely
combinatorial fashion. Let $S$ be an oriented smooth one-dimensional manifold 
without boundary and $f:S\to \R^2$ an immersion.

A double point $d\in\R^2$ of $f$ or the image in $\R^2$ of a critical
point of the composition $\begin{CD}S@>f>>\R^2@>{p_y}>>\R\end{CD}$
is called a {\it
critical point\/} of $f(S)$. A line passing through a critical point
of $f(S)$ and parallel to $x$-axis is called a {\it critical level.\/}

Assume that the immersion $f$ is {\it generic\/} in the sense that
\begin{enumerate}
\item it has only transversal double self-intersections,
\item its composition with the projection to $y$-axis has only
non-degenerate critical points,
\item no critical point of its composition with the projection
to $y$-axis is a double point and
\item each of its critical levels contains only one critical point.
\end{enumerate}

Fix a  real number $\GD>1$. A triple $a<b<c$ of points on a line is
called {\it $\GD$-symmetric\/} if 
$\GD^{-1}<\frac{c-b}{b-a}<\GD$.
It is easy to see that any horizontal line,
which meets $f(S)$ sufficiently close to a critical point, intersects
$f(S)$ in three points, which are not $\GD$-symmetric.

A generic immersion $f$ is said to be {\it $\GD$-regular\/} if there
are neighborhoods of the critical levels such that any horizontal line
which intersects $f(S)$ in a non-$\GD$-symmetric triple of points lies
in one of the neighborhoods of a critical level and two of the three
points are close to the critical point.

It is clear that for any generic immersion $f$ there exists 
sufficiently large $\GD$ such that $f$ is  $\GD$-regular.

\begin{figure}
\centerline{\epsffile{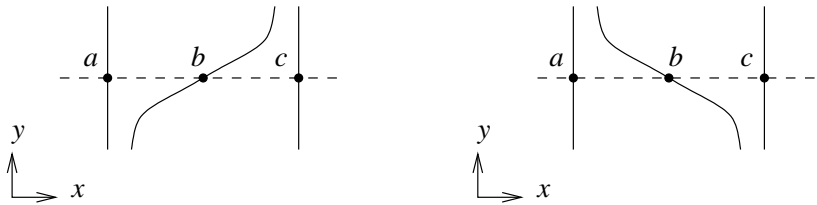}}
\caption{Appearance and disappearance of $\GD$-symmetric triples
in a slice of $\GD$-nonassociative immersion.}
\label{GD-tangle}
\end{figure}

A generic immersion $f$ with a finite number of critical levels is 
{\it $\GD$-nonassociative\/}
if the following conditions hold for any horizontal line containing a
$\GD$-symmetric triple $a<b<c$ of points of $f(S)$:
\begin{enumerate}
\item the line contains neither critical points nor other triples of
$\GD$-symmetric points of $f(S)$,
\item $(a,c)\cap f(S)=b$,
\item the $\GD$-symmetric
triple disappears in two different ways as the line moves up and down
with $\frac{c-b}{b-a}$ varying from $\GD^{-1}$ to $\GD$.
\end{enumerate}
The two possible types of a neighborhood of $[a,c]$ are shown in
Figure \ref{GD-tangle}.

Again, one can see that for any generic immersion $f$ there exists 
sufficiently large $\GD$ such that $f$ can be deformed by a diffeotopy
of the plane to a $\GD$-nonassociative immersion.

Non-associative immersions are related to nonassociative tangles
considered in \cite{Bar-Nat2}, \cite{Cartier}, which motivated our 
choice of this term.

The picture of a $\GD$-nonassociative immersion can be divided into
standard horizontal strips by lines separating the fragments 
containing $\GD$-symmetric triples and critical levels. 
Each of the strips should contain either a single critical point or 
a fragment shown in Figure \ref{GD-tangle}. 
This decomposition is referred to as a decomposition to elementary 
nonassociative fragments. 
It admits a purely combinatorial description in terms of bracketing, 
see \cite{Bar-Nat2}. 
An elementary fragment containing $\GD$-symmetric triple (i.e., a 
fragment shown in Figure \ref{GD-tangle}) is called an {\it
associator.\/}

\subsection{Elementary Characteristics of Regular and Nonassociative
Immersions}\label{s4.3}	Let $f:S\to\R^2$ be a generic immersion.
Denote by $M$ the number of maximum points of the composition of the
immersion and the projection to $y$-axis.

Denote by $X$ the number of the double points where 
either both branches are oriented upwards or downwards. If $S=\R^1$, 
this number is splitted as $X=X_++X_-$, where 
$X_+$ is the number of double points of $f$ where either both branches 
are oriented upwards and their intersection number is $+1$, or both 
branches are oriented downwards and the intersection number is $-1$.  

An immersion $f:\R^1\to\R^2$ such that $f(x)=(0,x)$ for $x\in \R$ with 
sufficiently large $|x|$ is called a {\em long curve}.  Thus, a long 
curve coincides at infinity with the standard parametrisation 
of the vertical axis.

\subsection{Casson Invariant via Nonassociative Diagram}\label{s4.4}
Consider a $\GD$-nonasso\-ci\-a\-tive long curve $\R\to\R^2$ and an 
associator $A$ appearing in its decomposition.  The three branches can 
be enumerated in two ways: from left to right and as their preimages 
appear in the source.  Denote by $\Gs(A)$ the element of $S_3$ which 
assigns the number of the branch counted according to the orientation 
on the source to the number of the same branch counted from left to 
right in the target.  Denote by $q(A)$ the number of branches of $A$ 
which are oriented upwards.  Define the sign $\Ge(A)$ of $A$ to be 
$(-1)^{q(A)}{\sign(\Gs(A))}$, if $A$ is as on the left hand side of 
Figure \ref{GD-tangle}, and $-(-1)^{q(A)}{\sign(\Gs(A))}$, if $A$ is 
as on the right hand side of Figure \ref{GD-tangle}.  Note that the 
sign $\sign(\Gs(A))$ and hence $\Ge(A)$ depend only of the cyclic 
order of the branches.  Thus $\Ge(A)$ is defined also for an 
associator in an immersion of $S^1$.  

For a nonassociative long curve and $\Go\in S_3$ put
$N(\Go)=\sum\Ge(A)$, where $A$ runs over all the
associators with $\Gs(A)=\Go$ appearing in the decomposition of the 
immersion. 

Let $N$ be the total algebraic number of 
associators. It is defined for a nonassociative immersion of either 
$\R^1$ or $S^1$. In the case of $\R^1$ it splits: $N=\sum_{\Go\in 
S_3}N(\Go)$.

\begin{thm}\label{4.4.A}
Let $K$ be a long knot whose projection to $xy$-plane is a
$\GD$-nonasso\-ci\-ative immersion for some $\GD>0$. Let $G$ be the
corresponding Gauss diagram. Then

\begin{equation}\label{4.4.A.2}
v_2(K)=\frac12\left\langle \vcenter{\epsffile{xleft.eps}}\ 
+\ \vcenter{\epsffile{xright.eps}}\ , G\right\rangle +\frac14(N(1)+N(1,3))
+\frac14X-\frac14M
\end{equation}
\begin{equation}\label{4.4.A.1}
v_2(K)=\frac12\left\langle \vcenter{\epsffile{xleft.eps}}\ 
+\ \vcenter{\epsffile{xright.eps}}\ , G\right\rangle +\frac14(N(2,3)+N(1,3,2))
+\frac12X_+
\end{equation}
\begin{equation}\label{4.4.A.3}
v_2(K)=\frac12\left\langle \vcenter{\epsffile{xleft.eps}}\ 
+\ \vcenter{\epsffile{xright.eps}}\ , G\right\rangle 
+\frac14(N(1,2)+N(1,2,3)) +\frac12X_-
\end{equation}
Here $X$, $X_+$, $X_-$ and $M$ are the characteristics of the 
projection of $K$ defined in Section \ref{s4.3}, and  $\left\langle
\vcenter{\epsffile{xleft.eps}} +\vcenter{\epsffile{xright.eps}}\ ,
G\right\rangle$ is the sum $\sum\Ge(c_1)\Ge(c_2)$  over all subdiagrams
of $G$ isomorphic to either $\vcenter{\epsffile{xleft.eps}}$, or
$\vcenter{\epsffile{xright.eps}}$, where $c_1$, $c_2$
are the chords of the subdiagram, see Section \ref{s1.2}. \end{thm}

Theorem \ref{4.4.A} is stated for a diagram of a {\em long} knot. Here 
is its reformulation for classical (closed) knots.

\begin{cor}\label{4.4.C} Let $K$ be a knot whose 
projection to $xy$-plane is a $\GD$-non\-as\-s\-o\-cia\-tive immersion 
for some $\GD>0$.  Let $G$ be the corresponding Gauss diagram.  Then 
\begin{equation}\label{4.4.C.eq}
v_2(K)=\frac14\left\langle \vcenter{\epsffile{x.eps}}\ , 
G\right\rangle +\frac1{24}N
+\frac18X-\frac1{24}M +\frac1{24}
\end{equation}
\end{cor}

\begin{proof} Obviously,
$$\left\langle \vcenter{\epsffile{xup.eps}}\ + \ 
\vcenter{\epsffile{xleft.eps}}\	+ \ \vcenter{\epsffile{xdown.eps}}\
+\ \vcenter{\epsffile{xright.eps}}\ , G\right\rangle=\left\langle 
\vcenter{\epsffile{x.eps}}\ , G\right\rangle$$
Now the result follows from \ref{T1}, \ref{Sym} and \ref{4.4.A}.
Appearance of $\frac1{24}$ is related to the fact that closing a
diagram long knot produces a new maximum point.
\end{proof}

\subsection{Casson Invariant via Regular Diagram}\label{s4.5}
Below by the index of a point $c$ with respect to a curve $\Gg$ we mean
the intersection number of $R_c$ and $\Gg$, where $R_c$ is
the open horizontal ray starting at $c$ and directed to the right.

Consider a $\GD$-regular long curve $f:\R\to\R^2$. The preimage of a
double point $d$ of $f$ divides $\R$ into three parts: two rays and a
segment. Denote by $i_{\text{int}}(d)$ the index of $d$ with respect to
the image of this segment under $f$.

\begin{figure}[hbt]
\centerline{\epsffile{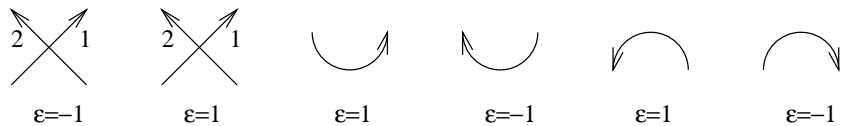}}
\caption{Signs of double and extremal points}
\label{signs}
\end{figure}

\begin{figure}[htb]
\centerline{\epsffile{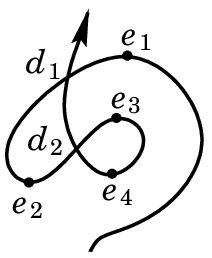}}
\caption{}
\label{fexample}
\end{figure}

Denote by $i_{\text{out}}(d)$ the
index of $d$ with respect to the image under $f$ of the union of the
rays.  Let $\Ge(d)$ be the
intersection number of the branches of $f(\R)$ at $d$. See Figure
\ref{signs}. In Figure
\ref{fexample}, $\Ge(d_1)=-1$, $\Ge(d_2)=1$, $i_{\text{int}}(d_1)=0$,
$i_{\text{out}}(d_1)=1$, $i_{\text{int}}(d_2)=-1$,
$i_{\text{out}}(d_2)=1$.

Put $I_{\text{out}}=\sum\Ge(d)i_{\text{out}}(d),$
$I_{\text{int}}=\sum\Ge(d)i_{\text{int}}(d),$
where the summations run over all double points $d$ of $f$.
For the curve in Figure \ref{fexample}, $I_{\text{int}}=-1$ and
$I_{\text{out}}=0$.

Local extrema of the composition of $f$ and 
the projection to the $y$-axis are called extremal points. At an
extremal point $e$ the curve $f$ goes either in a clockwise or a
counter-clockwise direction, see Figure \ref{signs}.
Let $\Ge(e)$ be $1$ in the counter-clockwise case and $-1$ otherwise.
In Figure \ref{fexample}, $\Ge(e_1)=\Ge(e_2)=1$ and
$\Ge(e_3)=\Ge(e_4)=-1$.

The preimage of an extremal point $e$ of $f$ divides $\R$ into two
rays. The curve $f(\R)$ is decomposed into the halves which are the
images of these rays.  Denote by $i_{\text{r}}(e)$ and
$i_{\text{l}}(e)$, respectively, the index of $e$ with respect to the
half of $f(\R)$  approaching $f(e)$ from the right and left,
respectively. Put $I_{\text{r}}=\sum\Ge(e)i_{\text{r}}(e),$
$I_{\text{l}}=\sum\Ge(e)i_{\text{l}}(e),$ where the summations run
over all extremal points $e$ of $f$. 
For the curve in Figure \ref{fexample}, $I_{\text{r}}=-1$ and
$I_{\text{l}}=0$.

\begin{thm}\label{4.5.A}
Let $K$ be a long knot whose projection to $xy$-plane and $y$-axis are
generic. Let $G$ be the
Gauss diagram corresponding to the projection to the
$xy$-plane. Then
\begin{equation}\label{4.5.A.2}
v_2(K)=\frac12\left\langle \vcenter{\epsffile{xleft.eps}}\ 
+\ \vcenter{\epsffile{xright.eps}}\ , G\right\rangle
-\frac14(I_{\text{out}}+I_{\text{r}})
+\frac14X-\frac14M
\end{equation}
\begin{equation}\label{4.5.A.1}
v_2(K)=\frac12\left\langle \vcenter{\epsffile{xleft.eps}}\ 
+\ \vcenter{\epsffile{xright.eps}}\ , G\right\rangle
+\frac12I_{\text{int}}
+\frac12X_+
\end{equation}
\begin{equation}\label{4.5.A.3}
v_2(K)=\frac12\left\langle \vcenter{\epsffile{xleft.eps}}\ 
+\ \vcenter{\epsffile{xright.eps}}\ , G\right\rangle 
-\frac14(I_{\text{out}}+I_{\text{l}}) +\frac12X_-
\end{equation}
Here $X$, $X_+$, $X_-$ and $M$ are the characteristics of the 
projection of $K$ defined in Section \ref{s4.3}, and  $\left\langle
\vcenter{\epsffile{xleft.eps}} +\vcenter{\epsffile{xright.eps}}\ ,
G\right\rangle$ is the sum $\sum\Ge(c_1)\Ge(c_2)$  over all subdiagrams
of $G$ isomorphic to either $\vcenter{\epsffile{xleft.eps}}$, or
$\vcenter{\epsffile{xright.eps}}$, where $c_1$, $c_2$
are the chords of the subdiagram, see Section \ref{s1.2}. \end{thm}

Theorem \ref{4.5.A} is stated for a diagram of a {\em long} knot.
However, it can be modified appropriately
giving rise to a formulation similar to Corollary \ref{4.4.C}.

For a generic immersion $f:S^1\to\R^2$ put $E=\sum \Ge(e)i(e)$, where
$e$ runs over extremal points of $f$ and $i(e)$ is the index of $e$
with respect to $f$. The preimage of a double point $d$ of $f$ divides
$S^1$ into two arcs. The curve $f(S^1)$ is decomposed into the halves
which are the images of these arcs.  One of them turns at $d$ in the
clockwise direction, the other one turns counter-clockwise. Denote by
$q_+(d)$ and $q_-(d)$, respectively, the index of $d$ with respect to
the former and latter, respectively. Put $Q=\sum(q_+(d)-q_-(d))$, where
$d$ runs over all double points of $f$.

\begin{cor}\label{4.5.B} Let $K$ be a knot whose 
projection to $xy$-plane and $y$-axis are generic.
Let $G$ be the  Gauss diagram corresponding to the projection to the
$xy$-plane.  Then 
\begin{equation}\label{4.5.C.eq}
v_2(K)=\frac14\left\langle \vcenter{\epsffile{x.eps}}\ , 
G\right\rangle -\frac1{24}E +\frac12Q
+\frac18X-\frac1{24}M +\frac1{24}
\end{equation}
\end{cor}

The proof is similar to the proof of \ref{4.4.C}.
Turn a closed curve to a long curve by cutting the leftmost string and
moving the cut points up and down. Clearly, $E$ turns into
$I_{\text{l}}+I_{\text{r}}$. Also, it is easy to check that $Q$ turns
$I_{\text{int}}-I_{\text{out}}$. Furthermore, $M$ increases by 1.
Now the result follows from \ref{T1}, \ref{Sym} and \ref{4.5.A}. \qed

\subsection{Digression: Relation to Arnold's Invariants of Plane
Curves}\label{s4.1.6}
Notice that in all the formulas of this Section there is a part
depending only on the knot projection. Moreover the rest of the
formula is common for all of the formulas:
$\left\langle2\ \vcenter{\epsffile{xup.eps}}+
2\ \vcenter{\epsffile{xleft.eps}}+
2\ \vcenter{\epsffile{xright.eps}}\ , G\right\rangle.$
Thus the parts of the formulas depending only on the plane curve
represent the same characteristic of the plane curve. Denote it by $I$.

It is easy to identify it with a linear combination $8St+4J^+$ of
invariants $St$ and $J^+$ of a generic immersion introduced by  Arnold
\cite{Arnold}.

Indeed, consider an ascending diagram of the unknot with the given
planar projection. Since $v_2$ of unknot is $0$,
$$I =-\left\langle2\
\vcenter{\epsffile{xup.eps}}+ 2\ \vcenter{\epsffile{xleft.eps}}+ 2\
\vcenter{\epsffile{xright.eps}}\ , G\right\rangle,$$
where $G$ is the corresponding Gauss diagram of the unknot.
The latter coincides with the Gauss diagram formula for
$4(2St+J^+)$ proved in \cite{Polyak}. 

$\frac14I$ coincides with the invariant which was extracted by Lin and
Wang \cite{LW} from Bar-Natan's integral formula \cite{Bar-Nat}
representing $v_2(K)$. 

All the formulas for $v_2(K)$ described above in this section can be
considered as formulas for $8St+4J^+$.

\section{New Configuration Spaces and Formulae}\label{s5}

\subsection{A Digression on the Degree of a Map}\label{s5.1}
The two classical methods for calculating the degree of a map 
discussed in Section \ref{s.straight} above, admit the following 
common generalization. Let $M$ and $N$ be  oriented smooth closed 
manifolds of dimension $n$.  Let $N$ be connected and $L$ be its
oriented smooth closed connected submanifold. Let $f:M\to N$ be a
differentiable map transversal to $L$. The orientations of $N$ and $L$
define an orientation of the normal bundle of $L$ in $N$. Because of
transversality, $f^{-1}(L)$ is a smooth submanifold of $M$. The normal
bundle of $f^{-1}(L)$ is naturally isomorphic to the pull back of the
normal bundle of $L$ and gets oriented. This orientation together with
the orientation of $M$ defines an orientation of $f^{-1}(L)$.
Therefore, both $L$ and $f^{-1}(L)$ are oriented smooth closed
manifolds of the same dimension and the map $g:f^{-1}(L)\to L$ defined
by $f$ has a well-defined degree. Obviously, this degree coincides 
with the degree of $f$.

It can be calculated by both of the classical methods. However the 
first method gives the expression for the degree literally coinciding 
with the one obtained by this method applied to the original map. 
The second method gives a new expression for the degree. In the case 
when $L$ is a point, this coincides with the expression obtained by 
the first method. 
In the case $L=N$, this coincides with the expression provided by the
second method applied to $f$. So, this is indeed a generalization of
both methods.

There are obvious generalizations of this observation. First, $M$
can be a stratified pseudomanifold. Then the transversality condition
is formulated as follows: the restrictions of $f$ to strata
of dimension $\ge \dim L$ are transversal to $L$. Second, a relative
situation can be considered: $M$ and $N$ are replaced by pairs $(M,M_0)$
and $(N,N_0)$ with $H_n(N,N_0)=\Z$ and $H_{\dim L}(L,L_0)=\Z$ where
 $L_0=L\cap N_0$.
Then the degree of a map $f:(M,M_0)\to(N,N_0)$ is equal to the degree
of the induced map
$$(f^{-1}(L),f^{-1}(L_0))\to(L,L_0).$$

Despite of apparent simplicity of this trick, it allows us to obtain 
several new geometrically interesting presentations of $v_2$.  We 
apply it to the map $\Gf:(\cal C,\cal S)\to (S^2\times S^2\times 
S^2,D)$ and various $L\subset S^2\times S^2\times S^2$.  However, this 
scheme is not easy to follow.  The first difficulty is related to the 
transversality condition.  The map has to be transversal to $L$ on 
each stratum of $\cal C$ (of all the dimensions).  The number of 
strata is rather large.  Moreover, in all the interesting cases there 
are strata on which the transversality condition is not satisfied.  
Another difficulty is that in the most interesting cases $L\subset D$.  
This problem is similar to the one we encountered in Section \ref{s4}.  
There we stepped back to a generic situation and then passed to limit.  
Here we can follow the same pattern, but prefer to consider the 
geometry related to $L$ in detail.  

The initial point of this consideration is still the structure of 
$\Gf^{-1}(L)$. However now we first take the intersection of 
$\Gf^{-1}(L)$ 
with the union $\cal C^0$ of all 6-dimensional strata of $\cal C$. 
Denote it by $\cal C^0_L$.
Then take the closure of $\cal C^0_L$ in $\cal C$. Denote the resulting 
space $\overline{\Gf^{-1}(L)\cap\cal C^0}$ by $\cal C_L$.  It is smaller 
and simpler than 
$\Gf^{-1}(L)$: even some high-dimensional strata of $\Gf^{-1}(L)$ do not 
show up.  This does not mean that we start over again from scratch.  We 
use the way how the high-dimensional strata of $\cal C$ are attached to 
each other.  Furthermore, since on these strata $\Gf$ is transversal to 
$L$, we can define the orientations using the scheme above.  Calculation 
of the degree is reduced to the case studied above via consideration of 
the preimage of a regular value.

\subsection{Locking the Free Point on a Chord}\label{s5.2}
Choose for $L$ the diagonal
$$\{(u_1,u_2,u_3)\in S^2\times S^2\times
S^2\mid u_2=u_3\}.$$					  

One can check that for a knot in general position the restriction of 
$\Gf$ to each 6-dimensional stratum of $\cal C$ is transversal to $L$. 
Let us identify $L$ with $S^2\times S^2$ by $(u_1,u_2,u_3)\mapsto (u_1,u_2)$
and denote by $\Gf_L:\cal C_L\to S^2\times S^2$ the map induced by 
$\Gf$. 

$\cal C_L$ is a four-dimensional pseudo-manifold. Its four-dimensional 
strata are the components of $\Gf^{-1}(L)\cap\cal C^0$. 
The components originated from $C_Y^0\times \Go$ can be identified with 
subspaces of $C^0_Y$ obtained by locking the free point $x_0$  on 
a line connecting $x_i$ and $x_j$ with $1\le i< j\le 3$. These subspaces 
are even closer to the initial motivation for introducing auxiliary 
strata, see Figure \ref{tripgraph}. 

\begin{figure}[htb]
\centerline{\epsffile{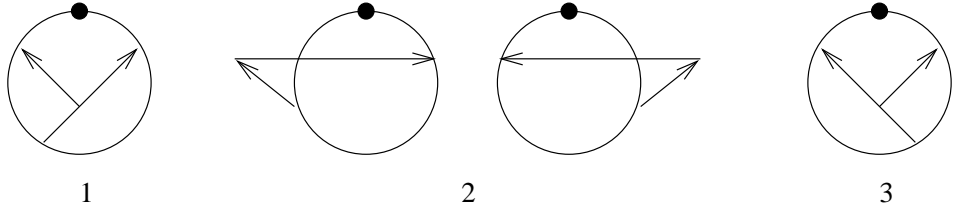}}
\caption{Diagrams representing $C^0_{Y,i}$, $i=1,2,3$.}
\label{fU}
\end{figure}

The strata related to $C_Y^0$ look as follows.
For $i=1,2,3$, denote by $C^0_{Y,i}$ the subspace of $C^0_Y$
defined by the condition \begin{enumerate}
\item $x_0$ lies on the line connecting $x_2$ and $x_3$ between
them, if $i=1$,
\item $x_0$ lies on the line connecting $x_1$ and $x_3$ outside
$[x_1,x_3]$, if $i=2$,
\item $x_0$ lies on the line connecting $x_2$ and $x_1$ between
them, if $i=3$.
\end{enumerate}
In the obvious sense, these spaces are associated with the diagrams
shown in Figure \ref{fU}.
Denote by $C_{Y,i}$ the closure of $C^0_{Y,i}$ in $C_Y$.

The intersection of $\Gf^{-1}(L)$ with $C_Y^0\times\Go$ can be identified
with:
\begin{enumerate}
\item $C_{Y,1}^0$ if $\Go=1$ or $(2,3)$;
\item $C_{Y,2}^0$ if $\Go=(1,2)$ or $(1,3,2);$
\item $C_{Y,3}^0$ if $\Go=(1,3)$ or $(1,2,3)$.
\end{enumerate}
Under these identifications, $\Gf_L$ turns into the maps extending the
ones defined by the following formulas on $C^0_{Y,i}$ with $i=1,2,3$
respectively:
\begin{enumerate}
\item $(x_1,x_2,x_3,x_0) \mapsto
\left(\frac{x_1-x_0}{|x_1-x_0|},\frac{x_3-x_0}{|x_3-x_0|}\right)$;
\item $(x_1,x_2,x_3,x_0) \mapsto
\left(\frac{x_0-x_2}{|x_0-x_2|},\frac{x_3-x_0}{|x_3-x_0|}\right)$;
\item $(x_1,x_2,x_3,x_0) \mapsto
\left(\frac{x_3-x_0}{|x_3-x_0|},\frac{x_1-x_0}{|x_1-x_0|}\right)$.
\end{enumerate}

Four other four-dimensional strata of $\cal C_L$ can be identified with
$C^0_X$. These strata are the intersections of
$\Gf^{-1}(L)$ with $C_X^0\times S^2\times \Go$, for $\Go=1,(2,3),(1,2),
(1,3,2)$.   On two of these strata $\Gf_L$ is identified with $\Gf^0_X$
and on two others, with $\Gf^0_X$ followed by the permutation of the
factors in $S^2\times S^2$. 

For the two remaining $\Go\in S_3$ (i.e., $(1,3)$ and
$(1,2,3)$), the intersections of $\Gf^{-1}(L)$ with $C^0_X\times
S^2\times \Go$ can be naturally identified with the product
$C^0_{II}\times S^2$, where 
\begin{equation}\label{C0II}C^0_{II}=\{(x_1,x_2,x_3,x_4)\in C^0_X\mid
x_1-x_3=\Gl(x_4-x_2) \text{ with } \Gl>0\}.\end{equation}
Under this 
identification,
$\Gf_L$ turns to $\Gf^0_{II}\times\id_{S^2}$, where
$$\Gf^0_{II}:C^0_{II}\to
S^2:(x_1,x_2,x_3,x_4)\mapsto\frac{x_1-x_3}{|x_1-x_3|}.$$

The identifications which have been made in the construction of $\cal C$ 
reduce the boundary of $\cal C^0_L$ in $\cal C_L$. The remaining 
3-dimensional strata of the boundary coincide with the 3-dimensional 
strata of $\cal C_L\cap \cal S$.
Denote the closure of these boundary strata by $\cal S_L$.

Observe that each 4-dimensional stratum of $\cal C_L$ described above 
appears twice, with the same mapping to $S^2\times S^2$.  This happens 
because $L$ is defined by $u_2=u_3$ and the permutation $(2,3)$, which 
acts in $\cal C$, commutes with $(2,3)$, which acts in $S^2\times 
S^2\times S^2$.  Therefore, we can quotient out $\cal C_L$ by the induced 
involution. The resulting space $\cal C_L/\Z_2$ has six 4-dimensional 
strata: $C^0_{Y,1}$, $C^0_{Y,2}$, $C^0_{Y,3}$, two copies of $C^0_X$,
and the stratum $S^2\times C^0_{II}$.    

Denote by $D^l\subset S^2\times S^2$ the set $$\{(u_1,u_2) \mid 
u_1=-u_2\}\cup \{(u_1,u_2)\mid u_1=\pm a\}\cup \{(u_1,u_2)\mid u_2=\pm 
a\}.$$ It is a bouquet of three copies of $S^2$.
Note that $D^l$ has codimension $2$ in $S^2\times S^2$, in contrast 
to $D$, which has codimension $1$ in $S^2\times S^2\times S^2$.  A 
straightforward modification of Theorem \ref{main_th} looks as follows. 

\begin{thm}\label{main_th_L} The space $\cal C_L/\Z_2$	
has a
well-defined fundamental class $[\cal C_L/\Z_2]\in H_4(\cal C_L/\Z_2,\cal 
S_L/\Z_2)$. The map
$\Gf_L:\cal C_L\to S^2\times S^2$ induces homomorphism
$$H_4(\cal C_L/\Z_2,\cal S_L/\Z_2)\to H_4(S^2\times S^2,D^l),$$
which maps $[\cal C_L/\Z_2]$ to $3v_2(K)[S^2\times S^2]$.\qed 
\end{thm}

\subsection{Configurations of Parallel Arrows}\label{} Aa part of
$\cal C_L/\Z_2$ coincides with $C^0_{II}\times S^2$, which is mapped 
to $S^2\times S^2$ by $\Gf_{II}\times\id_{S^2}$. 
This part of $\cal C_L/\Z_2$ differs in its nature 
from the rest of the main strata. It turns out that this part can be 
splitted out.

Observe that $C_{II}^0$ is
is the preimage of the diagonal under the map $\Gf^0:C^0_X\to S^2\times
S^2$.
The orientation of $C_{II}^0$ (defined by the orientation of $\cal C_L$)
can be computed as described in Section \ref{s5.1}. 
To a generic point $(x_1,x_2,x_3,x_4)\in C^0_{II}$, there correspond the
projection $\pi$ along $x_1-x_3$ and two crossings of this projection:
$c_{13}=\pi(x_1)=\pi(x_3)$ and $c_{42}=\pi(x_4)=\pi(x_2)$.
At such a point, one can take  as local coordinates $x_1$ and $x_3$.
The local degree of $\Gf_{II}^0$ with respect to the orientation defined
by the coordinate system $(x_1,x_3)$ is the sign $\Ge(c_{13})$ of $c_{13}$.
The local degree of $\Gf_{II}^0$ at this point is the product of this sign
and $\Ge(c_{42})$.
Hence the orientation of $C_{II}^0$ differs from the orientation defined
by the coordinate system $(x_1,x_3)$ by $\Ge(c_{42})$ .

Denote by $C_{II}$ the  closure of $C_{II}^0$ in $C_X$ and denote by
$\Gf_{II}$ the extension of $\Gf^0_{II}$ to $C_{II}$. The 1-strata of
$\p C_{II}$, whose images under $\Gf_{II}$ do not coincide with $a$, are
$\GS II_{i,i+1}=C_{II}\cap\GS X_{i,i+1}$,  $i=1,2,3$. There are embeddings
$$e_1:\GS II_{1,2}\to\GS II_{2,3}:(x_1,x_1,x_2,x_3)\to(x_1,x_2,x_2,x_3),$$ 
$$e_2:\GS II_{3,4}\to\GS II_{2,3}:(x_1,x_2,x_3,x_3)\to(x_1,x_2,x_2,x_3).$$
Let $\cal C_{II}$ be the quotient space of $C_{II}$ obtained by
identification of $\GS II_{2,3}$ with $\GS II_{1,2}\cup\GS II_{3,4}$
via these embeddings. It is easy to check that $\Gf_{II}$ defines a map
$\cal C_{II}\to S^2$.
Denote this map by the same symbol $\Gf_{II}$.

\begin{thm}\label{main_th_II}  $\cal C_{II}$ has a
well-defined fundamental class $[\cal C_{II}]\in
H_2(\cal C_{II},\p\cal C_{II})$. The map
$\Gf_{II}:\cal C_{II}\to S^2$ induces homomorphism
$H_2(\cal C_{II},\p\cal C_{II})\to H_2(S^2,a)$,
which maps $[\cal C_{II}]$ to $v_2(K)[S^2]$.
\end{thm}

\begin{proof} The local degree of $\Gf_{II}$ at a generic point
$(x_1,x_2,x_3,x_4)\in\cal C_{II}$ is $\Ge(c_{13})\Ge(c_{42})$.
Summing up the local degrees at all points of the preimage
of some regular value, we obtain  the right hand side of \eqref{v2}.
\end{proof}

Let us return to the space $\cal C_L/\Z_2$ of Section \ref{s5.2}.
Remove $C^0_{II}\times S^2$ from $\cal C_L/\Z_2$ and sue the 
edge of the cut by $e_1\times\id_{S^2}$ and $e_1\times\id_{S^2}$.
Denote the result by $\cal C^l$ (superindex $l$ stands for ``line''
alluding to locking the free point on a line).  Denote by $\Gf^l$ the
map $\cal C^l\to S^2\times S^2$ defined by $\Gf_L$.
Denote by $\cal S^l$ the subspace of $\cal C^l$ obtained from $\cal S_L$.
Combining Theorems \ref{main_th_L}	and \ref{main_th_II}, we get the
following theorem.

\begin{thm}\label{main_th_l} The space $\cal C^l$ has a
well-defined fundamental class $[\cal C^l]\in H_4(\cal C^l,\cal S^l)$.
The map
$\Gf^l:\cal C^l\to S^2\times S^2$ induces homomorphism $H_4(\cal
C^l,\cal S^l)\to H_4(S^2\times S^2,D^l)$, which maps $[\cal C^l]$
to $2v_2(K)[S^2]$.\qed  \end{thm}

\subsection{Application: New Integral Formula}\label{s5.2'}Theorem 
\ref{main_th_l}	gives rise to an integral formula for $v_2$ by choosing 
the standard volume form on $S^2\times S^2$, pulling it back by $\Gf^{l*}$ 
and integrating over $\cal C^l$. Put 
$$\Go=\frac{x\,
dy\wedge dz+y\,dz\wedge dx+z\,dx\wedge 
dy}{(x^2+y^2+z^2)^{3/2}}.$$ 
Recall that for a knot $K$ we denote by 
$C^0_{V}$ the space of 3-tuples $(x_1,x_2,x_3)\in K^3$ of
points ordered in the natural way defined by the orientation of $K$
with $*\ne x_1\ne x_2\ne x_3\ne *$. 

\begin{cor}\label{int-ch-formula}
\begin{multline}
v_2(K)=\int_{C_X^0}\Go(x_1-x_3)\wedge\Go(x_4-x_2) + \\
\frac12\int_{C_V^0}\int_{t\in(0,1)}
\Go(x_1-x_2)\wedge\Go((x_2-x_3)+t(x_1-x_2))+\\
\frac12\int_{C_V^0}\int_{t\in (0,1)}
\Go(x_2-x_3)\wedge\Go((x_1-x_2)+t(x_2-x_3))+\\
\frac12\int_{C_V^0}\int_{t\in(-\infty,0)\cup (1,\infty)}
\Go(x_3-x_1)\wedge\Go((x_1-x_2)+t(x_3-x_1)) 
\end{multline} \qed
 \end{cor}

\subsection{Slicing by Parallel Planes}\label{s5.3}	Now choose for $L$ 
the 3-torus $$\{(u_1,u_2,u_3)\in S^2\times S^2\times S^2\mid 
u_1\cdot a=u_2\cdot a=u_3\cdot a=0\}.$$
For a knot in general position the restriction 
of $\Gf$ to each 6-dimensional stratum of $\cal C$ is transversal to 
$L$. 
Denote by $\Gf_L:\cal C_L\to S^1\times S^1\times S^1$ the map induced by 
$\Gf$. Now $\cal C_L$ is a 3-dimensional pseudo-manifold, which has
twelve 3-dimensional strata. 

\begin{figure}[htb]
\centerline{\epsffile{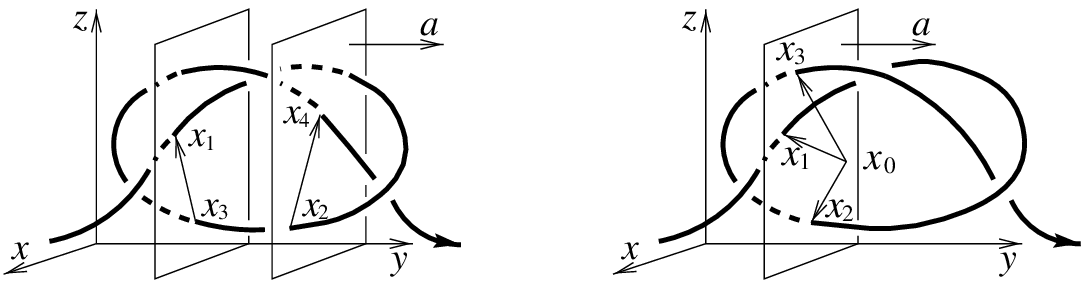}}
\caption{}
\label{slices}
\end{figure}

Six of them, which are 
$\Gf^{-1}(L)\cap\left(C^0_X\times S^2 \times \Go\right)$ with $\Go\in 
S_3$, can be identified with $C^a_X\times S^1\times \Go$, where 
$$C_X^a=\{(x_1,x_2,x_3,x_4)\in C_X^0\mid (x_1-x_3)\cdot 
a=(x_4-x_2)\cdot a=0\}.$$
In other words, $C_X^a$ is the space whose points are 4-tuples
$(x_1,x_2,x_3,x_4)$ of
points on the knot such that each of the pairs $(x_1,x_3)$ and $(x_2,x_4)$
lies in a plane orthogonal to $a$, see Figure \ref{slices}.

The other six strata, which are 
$\Gf^{-1}(L)\cap\left(C_Y\times \Go\right)$ with $\Go\in 
S_3$, can be identified with $C^a_Y\times \Go$, where 
$$C_Y^a=\{(x_1,x_2,x_3,x_0)\in C_Y^0\mid (x_1-x_0)\cdot 
a=(x_2-x_0)\cdot a=(x_3-x_0)\cdot a=0\}.$$
The space $C_Y^a$ consists of 4-tuples $(x_1,x_2,x_3,x_0)$ of points
lying in the same plane orthogonal to $a$ such that $x_1$, $x_2$ and
$x_3$ belong to $K$. See Figure \ref{slices}.

As above in Section \ref{s5.2}, 2-dimensional strata of the boundary
of $\cal C_L$ coincide with the 2-dimensional 
strata of $\cal C_L\cap \cal S$. They are contained in the closure of
 $\cup_{\Go} C^a_Y\times \Go$. A point of these strata can be
interpreted as a triple of points  contained in the same plane
orthogonal to $a$, which is tangent to the knot at one of these points
and intersects the knot in one of the other two points. Such a
configuration appears as a limit of 4-tuples belonging to $C_Y^a$ when
two of the points on the knot collide.

Denote the closure of these boundary strata by $\cal S_L$. 

Denote by $D^a\subset S^1\times S^1\times S^1$ the set
$$\{(u_1,u_2,u_3) \mid  u_1=-u_2\}\cup \{(u_1,u_2,u_3)\mid u_2=-u_3\}\cup
 \{(u_1,u_2,u_3)\mid u_3=-u_1\}.$$

Applying now the scheme described in Section \ref{s5.1}, we obtain the
following theorem.

\begin{thm}\label{main_th_a} $\cal C_L$ has a well-defined 
fundamental class $[\cal C_L]\in H_4(\cal C_L,\cal S_L)$. The map 
$\Gf_L:\cal C_L\to S^1\times S^1\times S^1$ induces homomorphism
$H_3(\cal C_L,\cal S_L)\to H_3(S^1\times S^1\times S^1,D^a)$,
which maps $[\cal C_L]$ to $6v_2(K)[S^1\times S^1\times S^1]$.\qed 
\end{thm}

\subsection{Yet Another Integral Formula}\label{s5.4} Theorem 
\ref{main_th_a}	gives rise to an integral formula for $v_2$ by choosing 
the standard volume form on $S^1\times S^1\times S^1$, pulling it back
by $\Gf_{L}^*$ 
and integrating over $\cal C_L$. Put 
$$\Ga=\frac1{2\pi}\frac{x\,
dz-zd\,x}
{x^2+z^2}.$$
This is the normalized length form on the circle in the $xz$-plane.

For fixed $t_1,t_2\in\R$, a 4-tuple $P=(w_1,w_2,w_3,w_4)$ of points
$w_i$ on $xz$-plane is called a {\it $(t_1,t_2)$-slice\/} if there
exists $(x_1,x_2,x_3,x_4)\in C_X^a$ such that $x_1$ and $x_3$ have
$y$-coordinates equal to $t_1$, $x_2$ and $x_4$ have $y$-coordinates
equal to $t_2$, and the projection of $x_i$ to $xz$-plane is $w_i$.
Denote by $p$ the number of $x_i$'s such that the projection to
$y$-axis of a positively oriented tangent vector to the knot  is
negative. For a fixed $t\in\R$, a triple $Q=(w_1,w_2,w_3)$ of points
$w_i$ on $xz$-plane is called a {\it $t$-slice\/} if there exists
$(x_0,x_1,x_2,x_3)\in C_Y^a$ such that $x_i$ has $y$-coordinate equal
to $t$ and the projection to $xz$-plane is $w_i$. Denote by $q$ the
number of $x_i$'s with $i=1,2,3$ such that the projection to
$y$-axis of a positively oriented tangent vector to the knot  is
negative.

\begin{cor}\label{liyboj}
\begin{multline}
v_2(K)=\int_{-\infty<t_1<t_2<\infty}\sum_{(t_1,t_2)-\text{slice}}
(-1)^p\Ga(w_1-w_3)\wedge\Ga(w_4-w_2) + \\
\int_{-\infty<t<\infty}\sum_{t-\text{slice}}
\int_{w_0\in\R^2\sminus\cup_{i=1}^3}
\Ga(w_1-w_0)\wedge \Ga(w_0-w_2) \wedge \Ga(w_3-w_0)
\end{multline} \qed
 \end{cor}

\begin{rem}The first integral is similar to the integral in the Kontsevich
formula for $v_2(K)$, see \cite{Konts} and \cite{Bar-Nat}. However, in our
formula we have instead of $dw/w$ only the imaginary part of $dw/w$.
Thus the contribution of the real part of $dw/w$ to the Kontsevich formula
is replaced by the second integral. A similar relation between the
associators in non-associative tangles and the contribution of $C_Y$
already appeared in Section \ref{s4.2}--\ref{s4.4}. This may shed light
onto yet non-understood relation between the integral formulas of
Kontsevich \cite{Konts}, \cite{Bar-Nat} and Bott-Taubes \cite{BT}.
\end{rem}


\begin{thebibliography}{99}


\bibitem{AM} S.~Akbulut and J.~McCarthy, {\em Casson's invariant for
oriented homology 3-spheres --- an exposition\/}, Princeton Math. Notes
36, Princeton University Press, 1990.
\bibitem{Arnold} V.I.Arnold, {\em
 Topological invariants of plane curves and caustics},
 University lecture series (Providence RI) {\bf 5}
 (1994)
\bibitem{AS}
S.Axelrod and I.M.Singer, {\em Chern-Simons Perturbatin Theory,} in
Proceedings of the XXth DGM Conference, edited by S.Catto and A. Rocha
(World Scientific, Singapore, 1992), 3--45; {\em Chern-Simons
Perturbation, II}, J. Diff. Geom. {\bf 39} (1994) 173--213.
\bibitem{Bar-Nat}
D.~Bar-Natan, {\em On the Vassiliev knot invariants\/}, Topology 
{\bf 34} (1995) 423--472.
\bibitem{Bar-Nat2} 
D.~Bar-Natan, {\em Non-associative tangles,} in {\em Geometric 
topology} (proceedings of the Georgia international topology 
conference), (W.~H.~Kazez, ed.), 139--183, AMS and International 
Press, Providence, 1997.
\bibitem{BT}
R.~Bott, C.~Taubes, {\em On the self-linking of knots}, J. Math. Phys.
{\bf 35} (1994), 5247--5287.
\bibitem{Cartier}
P.~Cartier, {\em Construction combinatoire des invariants de
Vassiliev-Kontsevich des noeuds}, C. R. Acad. Sci. Paris {\bf 316}
(1993), 1205--1210.
\bibitem{FM}
W. Fulton, R. MacPherson, {\em A Compactification of Configuration
Spaces}, Ann. Math. {\bf 139} (1994), 183--225.
\bibitem{Gilmer}
P.~Gilmer, {\em A method for computing the Arf invariants for links},
Quantum Topology, Series on Knots and Everything Vol. 3,
ed. L.~Kauffman and R.~Baadhio, World Sci., Singapore, 1993, 174--181.
\bibitem{GPV}
M. Goussarov, M. Polyak, O. Viro, {\em Finite Type Invariants of
Classical and Virtual Knots}, preprint math.GT/9810073.
\bibitem{Kauff}
Louis H. Kauffman, {\em On Knots}, Annals of Math. Studies 115,
Princeton University Press (1987).
\bibitem{Konts}
M.~Kontsevich,
{\em Vassiliev's knot invariants},
Adv. Sov. Math.
{\bf 16} (1993), 137--150.
\bibitem{Lannes}
Jean Lannes, {\em Sur les invariants de Vassiliev de degr\'e
infe\'rieur ou \'egal \`a 3},
 Enseign. Math. (2) {\bf 39} (1993) no. 3-4, 295--316.
\bibitem{Lannes'}
Jean Lannes, {\em Sur l'invariant de Kervaire des noueds classiques},
Comment. Math Helvetici, {\bf 60} (1985) 179--192.
\bibitem{LW}
Xiao-Song Lin and Zhenghan Wang, {\em Integral geometry of
plane curves and knot invariants},
J. Diff. Geom. {\bf 44} (1996), no. 1, 74--95.
\bibitem{Ng}
Ka Yi Ng, {\em Groups of ribbon knots},  Topology {\bf 37} (1998) 
441--458.
\bibitem{Piun}
S.~Piunikhin, {\em Combinatorial expressions for universal Vassiliev 
link invariant}, Comm. Math. Phys. {\bf 168-1} 1--22.
\bibitem{Poi} Sylvain Poirier, {\em Rationality results for the
configuration space integral of knots}, preprint math.GT/9901028.
\bibitem{PV}
Michael Polyak and Oleg Viro, {\em Gauss Diagram Formulas for
Vassiliev invariants}, Int. Math. Research Notices,
{\bf 11} (1994) 445--453.
\bibitem{Polyak}
Michael Polyak, {\em Invariants of curves and fronts via Gauss 
diagrams}, Topology {\bf 37} (1998) 989--1009.
\bibitem{P?} Michael Polyak, {\em On the algebra of arrow diagrams},
preprint http://www.math.tau.ac.il/\~{} polyak/.
\bibitem{Rolfsen} Dale Rolfsen, {\em Knots and links}, Publish or
Perish, Inc., 1990.
\bibitem{Thurston}
D.~Thurston, {\em Integral expressions for the Vassiliev knot 
invariants}, Harvard University senior thesis, April 1995.
\bibitem{V} 
V.A.Vassiliev, {\em Cohomology of knot spaces}, Theory of 
singularities and its applications (Providence) (V.I.Arnold, ed.) 
AMS, Providence, 1990.
\end{thebibliography}
\end{document}